\documentclass{amsart}
\usepackage{enumerate}

\theoremstyle{definition}
\newtheorem{ex}{Example}
\theoremstyle{plain}
\newtheorem{defin}{Definition}
\newtheorem{thm}[defin]{Theorem}
\newtheorem{lem}[defin]{Lemma}
\newtheorem{prop}[defin]{Proposition}
\newtheorem{cor}[defin]{Corollary}

\newcommand{\fps}[1]{\C [\! [#1 ]\! ]}
\newcommand{\fpsp}{\C [\! [Z,\zeta ]\! ]}
\newcommand{\CN}{\mathbb{C}^{N}}
\newcommand{\CNp}{\mathbb{C}^{N^\prime}}
\newcommand{\Cn}{\mathbb{C}^{n}}
\newcommand{\Cnp}{\mathbb{C}^{n^\prime}}
\newcommand{\Cd}{\mathbb{C}^{d}}
\newcommand{\Cdp}{\mathbb{C}^{d^\prime}}
\newcommand{\C}{\mathbb{C}}
\newcommand{\crb}[1]{\mathcal{V} (#1)}
\newcommand{\fcrvf}[1]{\mathcal{D}^{0,1}_{#1}}
\newcommand{\facrvf}[1]{\mathcal{D}^{1,0}_{#1}}
\newcommand{\fvf}[1]{\mathcal{D}_{#1}}
\newcommand{\Nn}{\mathbb{N}^{n}}
\newcommand{\N}{\mathbb{N}}
\newcommand{\dop}[1]{\frac{\partial}{\partial #1}}
\newcommand{\vardop}[2]{\frac{\partial#1}{\partial
#2}}

\newcommand{\Zp}{{Z^\prime}}

\newcommand{\cntocnp}{{(\CN, p_0)\to (\CNp,
p_0^\prime)}}
\newcommand{\innprod}[2]{\langle #1 , #2 \rangle}

\DeclareMathOperator{\spanc}{span _\C}
\DeclareMathOperator{\spanov}{span }
\DeclareMathOperator{\imag}{Im}

\DeclareMathOperator{\rk}{rk}

\begin{document}
\title[Holomorphic maps of real
submanifolds]{Holomorphic maps of real submanifolds
in Complex Spaces of different dimensions}
\date{November 9, 1999}
\author{Bernhard Lamel}
\address{Department of Mathematics, 0112\\
University of California, San Diego\\
La Jolla, CA 92092--0112}
\email{blamel@euclid.ucsd.edu}
\begin{abstract} We define a new local invariant
(called {\em degeneracy}) associated to  a triple
$(M,M^\prime,H)$, where
$M\subset\CN$ and
$M^\prime\subset\CNp$ are real submanifolds of
$\CN$ and $\CNp$, respectively, and $H: M\to
M^\prime$ is either a holomorphic map,
 a formal holomorphic map, or a smooth CR-map.
We use this invariant to find  sufficient
conditions under which finite jet dependence,
convergence and algebraicity results hold. 
\end{abstract}
\maketitle


\section{Introduction and statement of results}
In this paper, we discuss mappings of generic real
submanifolds in complex spaces of different
dimensions. We address the following specific
problems:
\begin{itemize}
\item Give conditions which ensure that a mapping
depends on its finite jet.
\item Give conditions under which a formal mapping
between real-analytic generic submanifolds is
convergent.
\item Give conditions under which
a map between algebraic
submanifolds is algebraic.
\end{itemize}

The first two questions have attracted considerable
attention in the equidimensional case, and quite
complete results have been obtained for the class
of finitely nondegenerate manifolds
(see \cite{BER3}, \cite{BER4}), and more recently,
for target manifolds of finite type in the sense
of D'Angelo (\cite{DAfintyp}) in \cite{BER5}.
Whether similar results hold for mappings of
generic submanifolds of spaces of different
dimension  is an intriguing problem which
leads to some new geometric notions. The
third question has also been answered in
terms of characterizing the algebraic
manifolds on which every holomorphic map is
algebraic (see especially \cite{Z2} and
\cite{BER2},
\cite{BHR1}, \cite{BR2}, \cite{H3}, \cite{W4},
\cite{M1}, \cite{SU1}, \cite{SS1}). We give a new
sufficient condition, which can be checked using
finitely many derivatives. 

 Our starting point is
the notion of
$(k_0,s)$-degeneracy.
This is a local invariant associated to the triple
$(M, M^\prime, H)$, where
$M\subset\CN$ and $M^\prime\subset\CNp$ are
generic $C^\infty$-submanifolds of $\CN$ and
$\CNp$, respectively, through 0,
and
$H:M\to M^\prime$ is a map (for example,
$C^\infty$-CR) which in loose terms measures
how ``flat'' $H(M)$ is as a
submanifold of
$M^\prime\subset\CNp$. 
The numbers $k_0$ and $s$ can be
defined (at $0$) as follows: If
$\rho^\prime_1,\dots \rho^\prime_{d^\prime}$ are
defining functions for $M^\prime$, $L_1,\dots,
L_n$ is a local basis for the CR-vector fields on
$M$,  and $H(0)=0$, then 
\begin{equation}\label{E:inform}
N^\prime - s = \max_k \dim_\C \spanc\Big\{ L^\alpha
\rho^\prime_{j,Z^\prime} (H(Z), \overline{H(Z)})
\big|_{Z=0}
\colon |\alpha|\leq k, 1\leq j \leq d^\prime \Big\},
\end{equation}
where for a multiindex $\alpha\in\Nn$ we write
$L^\alpha = L_1^{\alpha_1}\cdots L_n^{\alpha_n}$,
and $k_0$ is the least integer $k$ for
which the maximum dimension on
the right hand side of
\eqref{E:inform} is realized.
Here we write $N=n+d$, where $d$ is the
codimension of $M$. An extension of this definition
is given in   Section~\ref{section:formalities}
in the context of formal
submanifolds and formal maps,
which allows us a unified
treatment of real-analytic and
$C^\infty$-manifolds. This new notion is related
to the notion of {\em finite nondegeneracy} of a
real submanifold (first introduced in \cite{BHR1}),
and we explore this relationship further in
\ref{S:finitenondegeneracy}. 

Particularly satisfying is
the situation for  mappings
for which $s = 0$. We call
such mappings ``nondegenerate'',
or more specifically,
$k_0$-nondegenerate. These maps
fulfill a sufficient condition to
give a positive answer to all
three points above. The other maps
allowing for a further treatment
are the ones which are of
constant degeneracy (to be
defined in
Section~\ref{section:formalities}
as well).

 Let us
recall that a formal holomorphic map $H = (H_1,\dots ,
H_{N^\prime}): \CN\to\CNp$ at a point $p_0$ is  an
$N^\prime$-tuple of formal power series $H_j
(Z) =\sum_\alpha c^j_\alpha (Z - p_0)^\alpha$, and if
$H(p_0) = p_0^\prime
\in\CNp$, we write $H:(\CN, p_0)\to (\CNp,
p_0^\prime)$ for such a map. If $p_0\in M$, $p_0^\prime\in
M^\prime$ then we say that
$H:(\CN, p_0)\to (\CNp,
p_0^\prime)$ maps
$M$ into
$M^\prime $  if the following property is satisfied: If 
$\rho^\prime = (\rho^\prime_1,\dots ,
\rho^\prime_{d^\prime} )$ is a defining function
of $M^\prime$ and  $\rho = (\rho_1,\dots ,
\rho_{d} )$ is a defining function
of $M$ (where $d$ and $d^\prime$ are the codimensions of
$M$ and $M^\prime$, respectively), then there is a
$d^\prime\times d$ matrix $A$ of formal power series
such that $\rho^\prime (H(Z),\bar H (\zeta) ) = A
(Z,
\zeta ) \rho (Z,\zeta)$. (Here we
are  abusing notation: this
equation shall hold in the sense
of Taylor series.) 

Let  us recall that
we say that
$M$ is of finite type at $p$ (in the sense of
Kohn-Bloom-Graham) if the 
CR and the anti-CR vectors together with their
commutators of all lengths span the complexified
tangent space of $M$ at
$p$. We prove the following
theorems. If not stated explicitly
otherwise, all submanifolds are
assumed to be smooth and
connected.
\begin{thm}\label{T:nondegconv}
Let $M$, $M^\prime$ be generic real-analytic
submanifolds of $\CN$ and $\CNp$, respectively,
$p_0\in M$,  $M$ of finite type at $p_0$,
$p_0^\prime\in M^\prime$, and let
$H:\cntocnp$ be a formal holomorphic map which
maps $M$ into
$M^\prime$ and is
$k_0$-nondegenerate at $p_0$. Then there exists a
neighbourhood $U$ of $p_0$ in $\CN$ on which $H$ is
convergent. 
\end{thm}
For
the next theorem, we denote by
$j_{p_0}^k f$ the $k$-jet of $f$ at $p_0$.
\begin{thm}\label{T:jetdetnondeg1}
Let $M$, $M^\prime$ be generic real submanifolds of
$\CN$ and $\CNp$, respectively, $p_0\in M$, and $M$
of finite type at $p_0$. There exists an integer $K$
such that if
$H:M\to M^\prime$ and
$H^\prime: M\to M^\prime$ are $C^\infty$-CR
mappings which are 
$k_0$-nondegenerate at $p_0
\in M$ and $j_{p_0}^{Kk_0} H = j_{p_0}^{Kk_0}
H^\prime$, then $j_{p_0}^{l} H = j_{p_0}^{l}
H^\prime$ for all $l$. 
\end{thm}
\begin{thm}\label{T:jetdetnondeg2}
Let $M$, $M^\prime$ be generic real
submanifolds of $\CN$ and $\CNp$, respectively,
$p_0\in M$,  such that 
$M$ is of finite type at $p_0$. There exists an
integer $K$ such that if $H\colon U\to \CNp$ is a
holomorphic map defined on some neighbourhood $U$
of $p_0$ with $H(U\cap M)\subset
M^\prime$ and such that $H$ is
$k_0$-nondegenerate at
$p_0$, and $H^\prime$ is another holomorphic map
defined on some neighbourhood $U^\prime$ of $p_0$
with
$H^\prime (U^\prime\cap M)\subset M^\prime$ with 
\begin{equation}
\vardop{^\alpha H}{Z^\alpha} (p_0) = 
\vardop{^\alpha H^\prime}{Z^\alpha} (p_0), \qquad
|\alpha |\leq Kk_0,
\end{equation}
then $H = H^\prime$.
\end{thm}
Theorem~\ref{T:jetdetnondeg2} is an immediate
consequence of Theorem~\ref{T:jetdetnondeg1}. The
proof of Theorem~\ref{T:nondegconv} and
Theorem~\ref{T:jetdetnondeg1} is given in
Section~\ref{section:nondeg}.
Theorem~\ref{T:jetdetnondeg2} together with
the reflection principle in \cite{Boern2} yields
the following. 
\begin{cor}\label{C:jetdetnondeg}
Let $M$, $M^\prime$ be generic real-analytic
submanifolds of
$\CN$ and $\CNp$, respectively, $p_0\in M$, and $M$
of finite type at $p_0$. There exists an integer $K$
such that if
$H:M\to M^\prime$ and
$H^\prime: M\to M^\prime$ are $C^{Kk_0}$-CR
mappings which are 
$k_0$-nondegenerate at $p_0
\in M$ and $j_{p_0}^{Kk_0} H = j_{p_0}^{Kk_0}
H^\prime$, then both extend to holomorphic mappings
and
$H=H^\prime$. 
\end{cor}
Note that the notion of nondegeneracy makes
sense even for maps which are a priori only
smoothed to a finite order, so that the
statement of this corollary makes sense. The last
result we prove about nondegenerate maps is an
algebraicity theorem.
\begin{thm}\label{T:nondegalg} Let $M$ and
$M^\prime$ be algebraic generic submanifolds of
$\CN$ and
$\CNp$, respectively, $H$ a holomorphic map
defined on some connected neighbourhood $U$ of $M$ 
with $H(M)\subset M^\prime$ and
such that $H$ is 
$k_0$-nondegenerate at some point of $M$. Then
$H$ is algebraic. 
\end{thm}

Our next results are for hypersurfaces.  They are valid
either in the setting where
$N^\prime = N+1$, and the hypersurfaces are assumed to be
Levi-nondegenerate, or, where the
target hypersurface is strictly pseudoconvex  and
the source hypersurface is of finite type (and
there are no restrictions on
$N^\prime$). In the case of Levi-nondegenerate
hypersurfaces, we will consider  maps
$H$ which are (CR) transversal (the formal definition of
this property is given in
Definition~\ref{D:transversality}). We will refer
to the following properties in the theorems below:
\begin{enumerate}[(P1)]
\item $M^\prime$ is strictly
pseudoconvex at $p_0^\prime$.
\item $N^\prime = N+1$ and $M$ and $M^\prime$ are
Levi-nondegenerate at $p_0$ and $p_0^\prime$,
respectively, and $H$ is transversal at $p_0$.
\end{enumerate}

\begin{thm}\label{T:nponeconv}
Let $M$, $M^\prime$ be real-analytic
hypersurfaces in $\CN$ and $\CNp$, respectively,
$p_0\in M$, $p_0^\prime\in
M^\prime$, $M$ of finite type at $p_0$, and let
$H:\cntocnp$ be a formal holomorphic map
 of constant degeneracy which maps $M$ into
$M^\prime$. Then there exists a
neighbourhood $U$ of $p_0$ in $\CN$ on which $H$ is
convergent given that either (P1) or (P2)
holds. 

\end{thm}
Note that the case $N^\prime= N+1$ is very special, as the
following example shows. 
\begin{ex} Let $M\subset\CN$ be given by $\imag w =
\sum_{j=1}^{n} |z_j |^2$, and $M^\prime \subset \C^{N+2}$ be
given by $\imag w^\prime =  |z_{n+2} |^2 - 
\sum_{j=1}^{n+1} |z_j^\prime |^2 
$ (``adding a black hole''). Then the map
$(z_1,\dots,z_n,w)\mapsto (z_1 ,\dots,z_{n},f(z,w), f(z,w),w)$ maps $M$ into $M^\prime $
for {\em every} (formal) holomorphic map $f:\CN\to\C$.
\end{ex} 

This example also shows that in general, algebraicity
and dependence on jets of finite order for
Levi-nondegenerate hypersurfaces can only be
expected in the case
$N^\prime = N+1$ (without further restrictions on
the mappings, as for example nondegeneracy as
introduced above). 

\begin{thm}\label{T:jetdetnp1}
Let $M$, $M^\prime$ be real hypersurfaces in
 $\CN$ and $\CNp$, respectively, $p_0\in M$,
and $M$ of finite type at $p_0$. If
$H:M\to M^\prime$ and
$H^\prime: M\to M^\prime$ are $C^\infty$-CR
mappings which are  
constantly $(k_0,s)$-degenerate
at $p_0
\in M$ with $j_{p_0}^{2k_0} H = j_{p_0}^{2k_0}
H^\prime$, then $j_{p_0}^{l} H = j_{p_0}^{l}
H^\prime$ for all $l$, provided that either 
(P1) or (P2) holds.
\end{thm}
\begin{thm}\label{T:jetdetnp12}
Let $M$, $M^\prime$ be real hypersurfaces in
 $\CN$ and $\CNp$, respectively, $p_0\in M$,
and $M$ of finite type at $p_0$. If $H\colon U\to
\CNp$ is a holomorphic map defined on some
neighbourhood $U$ of $p_0$ with $H(U\cap
M)\subset M^\prime$ which is constantly
$(k_0,s)$-degenerate at
$p_0$, and $H^\prime$ is another holomorphic map
defined on some neighbourhood $U^\prime$ of $p_0$
with
$H^\prime (U^\prime\cap M)\subset M^\prime$ with 
\begin{equation}
\vardop{^\alpha H}{Z^\alpha} (p_0) = 
\vardop{^\alpha H^\prime}{Z^\alpha} (p_0), \qquad
|\alpha |\leq 2k_0,
\end{equation}
then $H = H^\prime$, provided that either (P1) or
(P2) holds.
\end{thm}
We also have the following algebraicity
result. Case (i) below is actually contained
in the results in \cite{Z2}.
\begin{thm}\label{T:np1alg} Let $M$ and
$M^\prime$ be algebraic hypersurfaces in
$\CN$ and
$\CNp$, respectively, $H$ a holomorphic map
defined on some connected neighbourhood $U$ of $M$
with $H(U\cap M)\subset M^\prime$.
Then
$H$ is algebraic, provided that either of
the following additional properties hold: 
\begin{enumerate}[(i)]
\item There exists a point $p_0$ in $M$ where $M$
is of finite type and $M^\prime$ is
strictly pseudoconvex at $H(p_0)$; 
\item $N^\prime = N+1$, and there exists a point
$p_0\in M$ at which $H$ is transversal, and
$M$ and
$M^\prime$ are Levi-nondegenerate at $p_0$ and
$H(p_0)$, respectively.
\end{enumerate} 
\end{thm}
Theorem~\ref{T:jetdetnp12} is again an immediate
consequence of Theorem~\ref{T:jetdetnp1}. In
the case $N^\prime = N+1$, with the assumptions of
the theorem,
$s=0$ or
$s=1$ (see Lemma~\ref{L:bound1}); the case $s=0$
is covered by Theorems~\ref{T:jetdetnondeg1} and
\ref{T:jetdetnondeg2}. The proofs of Theorems
\ref{T:nponeconv}, \ref{T:np1alg} and
 \ref{T:jetdetnp1} in the Levi-nondegenerate case
are given in Section
\ref{section:np1}. The proof for strictly
pseudoconvex hypersurfaces is given in
section~\ref{section:pcx}.

\section{Formal holomorphic maps of constant
degeneracy}\label{section:formalities}
\subsection{Some definitions}
\label{subs:formalmanifolds} In this section we want
to give a short review of some fundamental notions,
focusing on the parts which we shall need later on.
For a more thorough discussion of the notions
introduced here, see for example \cite{BER5}.

\subsubsection{Formal Submanifolds and Formal
 Maps}\label{subsub:formalthings} Consider the ring
of formal power series
$\fps{Z,\zeta}$ in the 2N indeterminates $(Z,\zeta
)=(Z_1,\dots, Z_N,\zeta_1,\dots , \zeta_N)$.
 A {\it formal generic submanifold}
$M\subset\CN$ of codimension $d$ is an ideal
$I\subset
\fps{Z,\zeta}$ which can be generated by $d$
formal power series  $\rho_1 (Z,\zeta)
,\dots,\rho_d (Z,\zeta )$ such that
$\bar
\rho_j (\zeta, Z)\in I$, $1\leq j \leq d$ and
$\rho_{1 Z} (0)\wedge \cdots \wedge \rho_{d Z}
(0)
\neq 0$ (where for $\phi\in\fps{Z,\zeta}$,
$\phi_Z= (\phi_{Z_1},\dots, \phi_{Z_N})$, and $\bar
\phi$ denotes the power series with barred
coefficients). A {\it formal holomorphic change of
coordinates}
$\tilde Z = F(Z)$ is the ring isomorphism between
$\fps{\tilde Z,\tilde
\zeta}$ and $\fps{Z,\zeta}$ induced by $\tilde
Z_j =  F_j(Z)$, $\tilde
\zeta_j= \bar F_j (\zeta)$, where $F=(F_1,\dots,
F_N)$ with $F_j\in \fps{Z}$ and these satisfy
\begin{equation}
\det
\left(\vardop{F_j}{Z_k}(0)
\right)_{\substack{1\leq j\leq N \\ 1\leq k \leq
N } } \neq 0.
\end{equation}

A {\it formal holomorphic map} $H:\CN\to\CNp$ is
an
$N^\prime$-tuple of formal power series in $Z$,
$H=(H_1,\dots, H_{N^\prime}) $ where
$H_j\in\fps{Z}$. This formal holomorphic map
induces a ring homomorphism $H^\sharp :
\fps{Z^\prime ,\zeta^\prime} \to \fps{
Z,\zeta}$ by $H^\sharp (Z_j^\prime) = H_j (Z)$,
$H^\sharp (\zeta_j^\prime)=\bar H_j (\zeta )$.
Given two formal submanifolds
$M\subset\CN$ and $M^\prime \subset\CNp$, of
codimension $d$ and $d^\prime$ respectively,
represented by their ideals
$I\subset\fps{Z,\zeta }$ and $I^\prime \subset
\fps{Z^\prime, \zeta^\prime}$ we
say that $H$ maps $M$ into $M^\prime$ (we will
be writing $H:M\to M^\prime$ to indicate that we
are in this situation) if
$H^\sharp (I^\prime) \subset I$. This is the case
if and only if for any generators $\rho =
(\rho_1,\dots,\rho_d)$ of $I$ and generators
$\rho^\prime =
(\rho_1^\prime,\dots,\rho_{d^\prime}^\prime)$ we
have that
\begin{equation}
\rho^\prime (H(Z),\bar H(\zeta )) = A(Z,\zeta)
\rho(Z,\zeta),
\end{equation}
where $A(Z,\zeta )$ is a $d^\prime\times d$
matrix of formal power series in $(Z,\zeta)$.
If in addition $d = d^\prime$ and $H$ is a
formal holomorphic change of coordinates, then
$\det A(0,0 ) \neq 0$, so that $A$ is an
invertible matrix of formal power series.

We work with formal power series since if we want
to handle $C^\infty$-submanifolds, then they
are a convenient way of keeping track of all
equations which we arrive from by repeated
differentiation. This brings us to the subject of
formal vector fields. A {\it formal vector field} 
$X$ is an operator of the form
\begin{equation}
X= \sum_{j=1}^N a_j (Z,\zeta ) \vardop{}{Z_j}
+\sum_{j=1}^N b_j(Z,\zeta ) \vardop{}{\zeta_j}.
\end{equation}
One checks that the formal
vector fields are exactly the derivations of the
$\C$-algebra
$\fps{Z,\zeta}$. $X$ is of {\it type} $(1,0)$ if
$b_j(Z,\zeta ) = 0$, $1\leq j\leq N$, and of
{\it type}
$(0,1)$ if $a_j(Z,\zeta ) = 0$, $1\leq j\leq N$, and
this distinction is invariant under formal
holomorphic changes of coordinates.  A formal
vector field is {\em tangent to M} if (as above,
$I$ denotes the ideal representing
$M$)
$XI\subset I$. 

A {\it formal CR-vector field} tangent to
$M$ is a formal vector field of type
$(0,1)$ tangent to
$M$. We write $\fcrvf{M}$ for the formal
CR-vector fields tangent to $M$. This is a
$\fps{Z,\zeta}$-module. Finally, we say that
$L_1,\dots,L_n$ (where $n=N-d$) is a {\em basis} of
the formal CR-vector fields tangent to $M$ if they
generate the quotient module $\fcrvf{M} /
I\fcrvf{M}$. For the anti-CR-vector fields tangent
to $M$ (that is, the formal vector fields of type
$(1,0)$ tangent to $M$) we write $\facrvf{M}$, and
set
$\fvf{M}=\fcrvf{M}\oplus \facrvf{M}$. For the
Lie algebra generated by $\fvf{M}$,  we write
$\mathfrak{g}$. We say that $M$ is of {\em finite
type} (at 0) if $\dim_\C \mathfrak{g}(o)= 2n+d$. 
\subsubsection{Formal normal coordinates}
\label{subsub:normalcoords} Let $M$ be a generic
formal submanifold of codimension $d$. Then after a
formal holomorphic change of coordinates we can
assume
$Z = (z,w) = (z_1,\dots,z_n,w_1,\dots , w_d ) \in
\Cn\times\Cd$ (for the corresponding $\zeta$ we
write $\zeta = (\chi,\tau)$) and that
$I$ is generated by $d$ functions $w_j - Q_j (z,
\chi,\tau)$, $j=1,\dots, d$, where
$Q_j\in\fps{z,\chi,\tau} $ fulfills
\begin{equation}\label{E:normal1}
Q_j (z,0,\tau) = Q_j (0, \chi,\tau ) = \tau_j,\quad
j=1,\dots,d.
\end{equation}
In that case, another useful set of generators is
given by  $\tau_j = \bar Q_j (\chi, z, w)$,
$j=1,\dots,d$. We shall call such coordinates {\em
formal normal coordinates}. We can use them in order
to parametrize $M$: Under this, we will understand
the ring isomorphism
\begin{equation}
\fps{z,w,\chi,\tau} / I \to \fps{z,\chi,\tau},
\quad  z\mapsto z, w\mapsto Q (z,\chi,\tau ),
\chi\mapsto \chi, \tau\mapsto \tau,
\end{equation}
or 
\begin{equation}\label{E:parambyqbar}
\fps{z,w,\chi,\tau} / I \to \fps{\chi,z,w},
\quad  z\mapsto z, w\mapsto w,
\chi\mapsto \chi, \tau\mapsto \bar Q
(\chi,z,w ).
\end{equation}
Note that a basis of CR-vector fields tangent to
$M$ is given by
\begin{equation}
L_k = \vardop{}{\chi_k} + \sum_{j=1}^{d} \bar
Q_{j\chi_k} (\chi, z, w) \vardop{}{\tau_j}, \quad
k=1,\dots, n.
\end{equation}
For a multiindex
$\alpha=(\alpha_1,\dots,\alpha_n)\in\Nn$ we write
$L^\alpha = L_1^{\alpha_1}\cdots L_n^{\alpha_n}$.
Let now $\phi(z,w,\chi,\tau) $ be a formal power
series. We want to relate the image of $\phi$ in the
parametrization of $M$ given by
\eqref{E:parambyqbar} with its derivatives along
CR-directions. Expand $\phi$ as a series in
$\fps{z,w} [\! [ \chi ]\! ]$:
\begin{equation}
\phi(z,w,\chi,\bar Q(\chi, z, w)) =
\sum_{\alpha\in\Nn} \frac{s_\alpha (z,w)}{\alpha
!}
\chi^\alpha.
\end{equation}
Now the $s_\alpha$ are obtained by partial
differentiation:
\begin{equation} \label{E:dchiasL}
s_\alpha (z,w) = \vardop{^{|\alpha |}}{\chi^\alpha}
\phi (z,w,\chi,\bar Q (\chi, z, w)) \Big|_{\chi = 0}
= L^\alpha \phi (z,w,0,w).
\end{equation}
The last equality is proved by induction on $|\alpha
|$ and its proof is left to the reader.

\subsubsection{Segre-mapppings and a finite type
criterion} Again, we are considering a generic
formal submanifold $M\subset\CN$ of codimension
$d$. Assume that formal normal coordinates $(z,w)$ as
in
\ref{subsub:normalcoords} have been chosen, along
with the corresponding (vector valued) function
$Q(z,\chi,\tau)\in\fps{z,\chi,\tau}^d$, fulfilling
\eqref{E:normal1}, such that $I$ (the ideal
associated to $M$) is generated by
$w_j - Q_j(z,\chi,\tau)$, $j=1,\dots,d$ (or,
equivalently, by
$\tau_j -\bar Q_j (\chi,z,w)$).  The {\it Segre
mappings} are the formal mappings $v^k:
\C^{kn}\to\CN$, for any integer $k$, defined by 
\begin{subequations}
\begin{align}\label{E:segremaps1}
&v^0 = (0,0), \qquad v^1(z) = (z,0),\\
&v^{2j}(z,\chi^1,\dots,z^{j-1},\chi^j)= \left(
z,Q(z,\chi^1,\bar Q(\chi^1,z^1,\dots, Q(z^{j-1},
\chi^j,0)\cdots ))\right), \label{E:segremaps2}\\
&v^{2j+1}(z,\chi^1,\dots ,\chi^j,z^{j})= \left(
z,Q(z,\chi^1,\bar Q(\chi^1,z^1,\dots, \bar Q(
\chi^j,z^{j},0)\cdots ))\right).\label{E:segremaps3}
\end{align}
\end{subequations}
We will write $(z,\chi^1,z^1,\dots)=(z,\xi)$,
where $\xi = (\chi^1,z^1,\dots)\in \C^{(k-1)n}$.
These mappings have the property that for
every
$k\geq 0$, and for every
$f\in I$,
\begin{equation}
\label{E:segremapiszero}
f(v^{k+1}(z,\xi), \bar v^k (\xi)) = 0.
\end{equation}
 We shall need the finite type criterion of
Baouendi, Ebenfelt and Rothschild, which we state
here for reference; see e.g. \cite{BER5} or
\cite{BERbook}.
 For a formal mapping $F$, $\rk{F}$  is the rank
of the Jacobian matrix of $F$ over the quotient
field of the ring of formal power series involved. 
\begin{thm}\label{T:bercrit} Let $M$ be a formal
generic submanifold of $\CN$ of codimension $d$.
Then, $M$ is of finite type at $0$  if and only if
there exists a $k_1\leq d+1$ such that $\rk (v^k) =
N$ for
$k\geq k_1$. Moreover, if $M$ is real-anlaytic and
of finite type at $0$, then there exists
$(z_0,\xi_0)\in
\C^n\times \C^{(2k_1 - 1)n}$ arbitrarily close to
the origin such that $v^{2k_1}(z_0,\xi_0 ) = 0$ and
the rank of the Jacobian matrix of $v^{2k_1}$ at
$(z_0,\xi_0 )$ is $N$. 
\end{thm}
Note the following consequence for real-analytic
submanifolds: There are
$(z_0, \xi_0 )\in \C^n\times \C^{2k_1}$, arbitrarily
close to
$0$, such that the function $v^{2k_1}$ has a
holomorphic right inverse $\psi:\CN\to \C^n\times
\C^{2k_1} $ defined in a neighbourhood of $0\in
\CN$, such that $\psi(0) = (z_0,\xi_0 )$ and
$v^{2k}(\psi(Z)) = Z$. This follows from the
inverse function theorem and Theorem~\ref{T:bercrit}.

\subsection{Constant-rank submodules}
We are now considering a free module of
rank $k$ over $\fps{Z,\zeta }$. Let
$E\subset \fps{Z,\zeta }^k $ be a submodule
and $I\subset \fps{Z, \zeta}$ an ideal.
\begin{defin}\label{D:constant rank}
The submodule $E\subset \fps{Z,\zeta }^k $ is
of constant rank $l= \dim_\C E(0)$ over $I$ if
for all $e_1,\dots,e_l,e_{l+1}\in E$, where $e_j =
(e_j^1,\dots,e_j^k)$, every subminor of length
$l+1$ of the matrix $(e_m^n)_{1\leq m \leq l+1,
1\leq n
\leq k}$ is an element of $I$.
\end{defin}

Definition~\ref{D:constant rank} is motivated by the
following observation. Consider a
$C^\infty$-submanifold $M\subset\CN$. By
taking the Taylor series of the defining
functions of $M$ at a point, one obtains a formal
submanifold represented by some ideal $I$. Consider
a vector bundle of rank
$l$ over $M$, embedded in $\C^k$. Its sections
will be a submodule of $\Gamma (M, \C^k)\cong
C^\infty (M)^k$ (the isomorphism is determined
by a choice of coordinates). Taking the Taylor
expansion of the components, we get a
submodule $E\subset \fps{Z,\zeta}^k$ which
fulfills the condition of
Definition~\ref{D:constant rank} over $I$. 

We may refer to a submodule of
constant rank $l$ as a formal vectorbundle over
$M$. This is also highlighted by the following
Lemma, which can be thought of as a characterization
of the bases of sections of a formal vectorbundle.
\begin{lem} \label{L:formalbundle} Suppose that
$E$ is a submodule of
$\fps{Z,\zeta}^k$, with $\dim_\C E(0) =l$. Then
$E$ is  of constant rank
$l$ over
$I\subset
\fps{Z,\zeta}$ if and only if $E$ has the following
property: If
$v_1,\dots, v_l\in E$ are vectors such that $v_1
(0),\dots v_l(0)$ form a basis of $E(0)$, then
$v_1,\dots , v_l$ generate $E$ up to vectors all of
whose components are elements of $I$.  
\end{lem}
\begin{proof} The ``if'' direction is trivial; so
assume that $E$ is of constant rank $l$. After
reordering, we can assume that if we write
$v_j = (v_j^1,\dots , v_j^k)$, the matrix
$(v_m^n)_{1\leq m,n \leq l}$ is invertible. Hence,
given any $e = (e^1,\dots, e^k)\in E$, we
can find $a_1,\dots, a_l\in\fps{Z,\zeta}
$ such that $\sum_{m=1}^l a_n v_m^j = e^j$,
$j=1,\dots, l$. Now consider $e^\prime = e- \sum_j
a_j v_j\in E$. We want to show that the components of
this vector are elements of $I$; this is clear for
the first $l$ components (which are 0, after all).
Taking a subminor of length
$l+1$ of the matrix
\begin{equation}
\begin{pmatrix}
v_1^1 &\dots & v_1^l &\dots &v_1^k \\
\vdots& & \vdots && \vdots \\
v_l^1 & \dots &v_l^l &\dots &v_l^k \\
{e^\prime}^1 &\dots &{e^\prime}^l&\dots
&{e^\prime}^k  
\end{pmatrix}
\end{equation}
which contains the first $l$ columns and developing
it along the last row, by assumption we have that
\begin{equation}\label{E:detinI}
\pm \begin{vmatrix}
v_1^1 &\dots & v_1^l \\
\vdots& & \vdots \\
v_l^1 & \dots &v_l^l
\end{vmatrix} {e^\prime}^j \in I, \quad l+1\leq
j\leq k
\end{equation}
which implies ${e^\prime}^j \in I$ since the
determinant in \eqref{E:detinI} is a unit in
$\fps{Z,\zeta}$.
\end{proof}

\subsection{Degeneracy of a Formal Holomorphic Map}
\label{sub:formaldegeneracy}
Let $H:M\to M^\prime$ be a formal holomorphic map
between the formal submanifolds $M\subset\CN$ and
$M^\prime\subset\CNp$. Choose a basis of formal
CR-vector fields
$L_1,\dots,L_n$ tangent to $M$ and generators
$\rho_1^\prime
(Z^\prime,\zeta^\prime ),\dots,\rho_{d^\prime}^\prime
(Z^\prime ,\zeta^\prime )$ of the ideal $I^\prime$
representing
$M^\prime$. The formal series
$\rho^\prime_1(H(Z),\bar H(\zeta ) ),\dots
,\rho^\prime_{d^\prime} (H(Z),\bar H(\zeta ) ) $ are
then elements of $I$ (since $H:M\to M^\prime$).
For $\rho^\prime_j(Z^\prime ,\zeta^\prime )$, we
denote by
$\rho_{j,Z^\prime}^\prime
(Z^\prime,\zeta^\prime)$ the {\em complex
gradient}
\begin{equation}
\rho_{j,Z^\prime}^\prime (Z^\prime,\zeta^\prime)
= \left( \vardop{\rho_{j}^\prime
(Z^\prime,\zeta^\prime)}{Z_1^\prime},\dots ,
\vardop{\rho_{j}^\prime
(Z^\prime,\zeta^\prime)}{Z_{N^\prime}^\prime}
\right)
\end{equation}
and we usually think of it as a {\em row} vector.
So $\rho_{j,Z^\prime}^\prime
(Z^\prime,\zeta^\prime)\in
\fps{Z^\prime,\zeta^\prime}^{N^\prime}$, and
$\rho_{j,Z^\prime}^\prime (H(Z),\bar H (\zeta ))
\in \fps{Z,\zeta}^{N^\prime}$. We define
an ascending chain of submodules
$E_k\subset\fps{Z,\zeta}^{N^\prime}$ by
\begin{multline}\label{E:theek}
E_k = \spanov_{\fps{Z,\zeta}} \{ L_1\cdots L_r
\rho_{j,Z^\prime}^\prime (H(Z),\bar H(\zeta
))\colon\\ L_1,\dots, L_r \text{
formal CR-vector fields tangent to } M, 0\leq
r\leq k, 1\leq j \leq d^\prime\}.
\end{multline}
The
chain $E_k(0)$ of subspaces of $\CNp$ will
stabilize, say at the index
$k_0$, that is, $E_k(0) =E_{k_0}(0)$ for $k\geq k_0$,
and
$E_{k_0 - 1}(0)
\neq E_{k_0}(0)$. So will the chain of submodules $E_k \subset
\fps{Z,\zeta}^{N^\prime}$, 
since $\fps{Z,\zeta}^{N^\prime}$ as a free module
over a Noetherian ring is Noetherian itself; say
$E=\cup_k E_k$, and there is some $k_0^\prime \geq
k_0$ such that $E_k = E$ for $k\geq k_0^\prime$.
\begin{defin}\label{D:formaldegeneracy}
With $E_{k_0}$ defined above, let $s=N^\prime -
\dim_\C E_{k_0} (0)$. We then say that $H$ has
formal degeneracy $(k_0, s)$ (at $0$) or shortly,
that
$H$ is
$(k_0, s)$-degenerate (at $0$). If $E$ is of
constant rank $N^\prime - s$ over $I$, we say that
$H$ has constant degeneracy $s$ and that $H$ is
constantly
$(k_0,s)$-degenerate. Furthermore, if $s=0$, we
say that $H$ is $k_0$-nondegenerate. Since in that
case, $E = \fps{Z,\zeta}^{N^\prime}$, $H$ has
automatically constant degeneracy $0$.
 \end{defin}
Note that if $H$ is of constant degeneracy, then the
submodule $E_{k_0}$ will actually generate $E$ up to
vectors whose components are in $I$ (this follows
from matrix manipulations, see Lemma
\ref{L:formalbundle}). However, in general, we do
not know whether $k_0 = k_0^\prime$. Clearly, the
degeneracy
$s$ fulfills the inequality
$0\leq s\leq N^\prime - d^\prime$.
Without further restrictions on
$M$,
$M^\prime$ and $H$, this is the best we can hope
for, as the example of a Levi-flat submanifold as
the target shows: If $M^\prime$ is defined by the
equations $\imag w_1^\prime = \cdots = \imag
w_d^\prime = 0$, where $Z^\prime =
(z^\prime,w^\prime) = (z_1^\prime,\dots ,
z_n^\prime,w_1^\prime,\dots ,w_d^\prime)$ are
coordinates in
$\CNp$, then every map $H$ which maps into this
manifold has constant degeneracy
$N^\prime - d^\prime$.

We will now show that
Definition~\ref{D:formaldegeneracy} is actually
independent of choices of formal coordinates and
generators. First consider a different set of
generators $\tilde \rho^\prime = (\tilde
\rho^\prime_1,\dots ,\tilde \rho^\prime_{d^\prime}
)$. Then there is an invertible $d^\prime \times
d^\prime$ matrix $A=(a_{jk})$ of formal power series
in
$\fpsp$ such that $\tilde \rho^\prime
(Z^\prime,\zeta^\prime)  = A(Z^\prime,\zeta^\prime)
\rho(Z^\prime,\zeta^\prime)$. Taking the complex
gradient, we obtain
\begin{equation}
\label{E:cplgrad} \tilde \rho^\prime_{j,Z^\prime}
 = \sum_{k=1}^{d^\prime}
\rho_k^\prime  a_{jk, \Zp}  +
\sum_{k=1}^{d^\prime} a_{jk} \rho^\prime_{k,
Z^\prime}, \quad 1\leq j \leq d^\prime.
\end{equation}
 Now the first sum in \eqref{E:cplgrad} is a vector
whose entries are elements of $I^\prime$. We write
$\phi_j (Z,\zeta) = \rho_{j, \Zp}^\prime
(H(Z),\bar H (\zeta ))$, $\phi =
(\phi_1,\dots,\phi_{d^\prime})$, use the same
notation for $\tilde \rho^\prime_{Z^\prime}$, and
set
$B(Z,\zeta ) = A(H(Z),\bar H (\zeta ))$; note that
$B$ is an invertible $d^\prime \times d^\prime$
matrix of formal power series in $\fps{Z,\zeta}$.
Then pulling
\eqref{E:cplgrad} to
$M$ via
$H$, we see that
\begin{equation}
\label{E:cplxgradonM}
\tilde \phi_j = v_j + B_j \phi, \quad 1\leq j
\leq d^\prime,
\end{equation}
where $v_j$ is a vector  in
$I\fps{Z,\zeta }^{N^\prime}$. \eqref{E:cplxgradonM}
implies that
$\tilde E_k = E_k$ modulo $I
\fps{Z,\zeta}^{N^\prime}$ for each $k$, which implies
that
$\dim
\tilde E_k(0) = \dim E_k(0)= e$ and that
$\tilde E_k$ is of constant rank e if and only if
$E_k$ is. This shows that the choice of defining
function does not matter.

Note that Definition~\ref{D:formaldegeneracy} is
independent of the choice of formal holomorphic
coordinates in
$\CN$; that follows easily from the fact that such a
biholomorphic change of coordinates $F$ pushes formal
CR-vector fields tangent to $M$ to formal
CR-vector fields tangent to $F(M)$. The
independence from choice of formal holomorphic
coordinates in $\CNp$ is  proved in the next
Lemma.
\begin{lem} Suppose $M\subset\CN$,
$M^\prime\subset\CNp$ are generic formal
submanifolds of $\CN$ and $\CNp$, respectively, and
that $H:M\to M^\prime$ is a formal holomorphic map.
When we change formal holomorphic coordinates in
$\CNp$ by $\tilde Z^\prime = F(Z^\prime)$, and denote
the space defined by \eqref{E:theek} in the $\tilde
Z^\prime$-variables by $\tilde E_k$, then\
\begin{equation}
\label{E:ektransform}
\tilde E_k \left( \vardop{F}{Z^\prime}
(H(Z))\right) = E_k.
\end{equation}
\end{lem}
\begin{proof} We choose generators
$\tilde\rho^\prime$ for the ideal
$\tilde I^\prime$ representing $M^\prime$ in the
coordinates $\tilde Z^\prime$. Then we can
choose
$\rho^\prime = F^\sharp
\tilde\rho^\prime$ as generators for
$I^\prime\subset\fpsp$. We now take the complex
gradient, and use the chain rule to obtain
\begin{equation}\label{E:gradtrans}
\rho_{j,Z^\prime}^\prime(Z^\prime ,\zeta^\prime ) =
\left(\tilde
\rho_{j} (F(Z^\prime),\bar F(\zeta^\prime)
\right)_{Z^\prime} = \tilde \rho_{j,\tilde
Z^\prime} (F(Z^\prime),\bar F(\zeta^\prime ))
\vardop{F}{Z^\prime} (Z^\prime).
\end{equation}
Pulling \eqref{E:gradtrans} to $M$ and applying
CR-vector fields tangent to $M$, we obtain
\eqref{E:ektransform}, since all (0,1)-vector
fields annihilate the entries of the matrix
$\vardop{F}{Z^\prime} (H(Z))$.
\end{proof}
Now the transformation
$v=(v_1,\dots,v_{N^\prime})\mapsto v A$, where $A$
is an invertible $N^\prime\times N^\prime$ matrix of
formal power series in
$(Z,\zeta)$, is a module isomorphism of
$\fps{Z,\zeta}^{N^\prime}$ which maps
$I\fps{Z,\zeta}^{N^\prime}$ into itself. So
Definition~\ref{D:formaldegeneracy} is in fact
independent of all the choices made there. The next
Lemma, the proof of which we leave to the reader,
gives a means of actually computing with
Definition~\ref{D:formaldegeneracy}.
\begin{lem}\label{L:deginbasis} Suppose $M\subset\CN$,
$M^\prime\subset\CNp$ are generic formal
submanifolds of $\CN$ and $\CNp$, respectively, and
that $H:M\to M^\prime$ is a formal holomorphic map.
Let $L_1,\dots, L_n$ be a basis of the CR-vector
fields tangent to $M$ (that is, a set of generators
of $\fcrvf{M}/ I\fcrvf{M}$,
see~\ref{subsub:formalthings}). For a multiindex
$\alpha\in\Nn$, let $L^\alpha = L_1^{\alpha_1}\cdots
L_n^{\alpha_n}$. Furthermore, let generators
$\rho^\prime = (\rho_1^\prime,\dots
\rho_{d^\prime}^\prime)$  of $I^\prime$ be chosen.
Define the submodules
\begin{equation}\label{E:deginbasis}
F_k = \spanov_{\fps{Z,\zeta}} \{ L^\alpha
\rho^\prime_{j,Z^\prime} (H(Z),\bar H (\zeta ))
\colon |\alpha |\leq k, 1\leq j\leq d^\prime\}
\end{equation}
of $\fps{Z,\zeta}^{N^\prime}$. Let $E_k$ be defined
by \eqref{E:theek}. Then $E_k(0) = F_k(0)$, and
$E_k$ is of constant rank over $I$ if and only if
$F_k$ is; hence, in order to determine the
degeneracy of $H$, it suffices to consider the
$F_k$.
\end{lem}
We now want to give a different characterization of the degeneracy
$s$ in the case of constant degeneracy. For that, we
will formulate the conditions of
Definition~\ref{D:formaldegeneracy} in terms of
formal normal coordinates (see
\ref{subsub:normalcoords}). So $I^\prime$ is
generated by the functions $\rho^\prime_j =
w_j^\prime - Q_j^\prime (z^\prime, \chi^\prime,
\tau^\prime)$,
$j = 1,\dots, d^\prime$, and we will write
$H=(f,g)=(f_1,\dots, f_{n^\prime} , g_1,\dots,
g_{d^\prime})$ for $H$ in these normal coordinates.
The complex gradient is easily computed to be
\begin{equation}\label{E:cplxgradnormal}
\rho^\prime_{j,Z^\prime} = (-Q^\prime_{jz_1},\dots, -
Q^\prime_{jz_{n^\prime}}, e_j), \quad 1\leq j
\leq d^\prime, 
\end{equation}
where $e_j$ is the $j$th unit vector in
$\C^{d^\prime}$. In particular, the last
$d^\prime$ entries of any CR-derivative of
length bigger than $0$ of any
$\rho_{j,Z^\prime}^\prime (H(Z),\bar H (\zeta))$ will be
$0$. We will write $Q_{j z_k^\prime}^\prime
(f(z,w),\bar f (\chi,\tau) , \bar g (\chi,\tau )) =
\phi_j^k (Z,\zeta)$. So $H$ is of constant
degeneracy $s$ if and only if after reordering the
$z$'s, we have that there exist $t = n^\prime - s$
multiindeces
$\alpha^1,\dots ,\alpha^t$, and
integers
$l_1,\dots, l_t$, $1\leq l_j \leq d^\prime$ such that
the vectors $(L^{\alpha^j} \phi_{l_j}^1,\dots,
L^{\alpha^j} \phi_{l_j}^t)$, $1\leq j \leq t$,
evaluated at
$0$, form a basis of $\C^t$, and for all
multiindeces
$\beta$, all $k$,
$t+1\leq k \leq n^\prime$, and all $l$, $1\leq l
\leq d^\prime$, the determinant
\begin{equation}\label{E:matrix1}
\begin{vmatrix}
L^{\alpha^1} \phi_{l_1}^1 &\dots &L^{\alpha^1} \phi_{l_1}^t
&L^{\alpha^1} \phi_{l_1}^k \\
\vdots &\ddots &\vdots &\vdots \\
L^{\alpha^t} \phi_{l_t}^1 &\dots &L^{\alpha^t} \phi_{l_t}^t
&L^{\alpha^t} \phi_{l_t}^k \\
L^{\beta} \,\, \phi_{l}^1 &\dots &L^{\beta} \,\,\phi_{l}^t
&L^{\beta}\,\, \phi_{l }^k
\end{vmatrix} \in I.
\end{equation}
More specifically, if $H$ is constantly
$(k_0,s)$-degenerate, one of the $\alpha_t$ must
have length $k_0$. We can use
this to formulate the following  technical result:
\begin{lem}\label{L:constdeg}
Assume that normal coordinates $(z,w)$ and
$(z^\prime, w^\prime)$ have been chosen for $M$ and
$M^\prime$, respectively, that $L_1,\dots,
L_{n}$ is a basis for the CR-vector fields
tangent to $M$, that $w_j^\prime -
Q_j^\prime(z^\prime ,\chi^\prime,\tau^\prime)$ are
generators of
$I^\prime$ as in \ref{subsub:normalcoords}, and
that
$H:M\to M^\prime$ is of constant degeneracy $s$. Let
$t=n^\prime -s$, and write $H^\sharp Q_{j
z_k^\prime}^\prime = \phi_j^k$. We can choose
$t$ multiindeces $\alpha^1,\dots \alpha^t$, and
integers $l_1,\dots,l_t$, such that after reordering
the
$z^\prime$'s  the vectors $(L^{\alpha^j}
\phi_{l_j}^1,\dots, L^{\alpha^j} \phi_{l_j}^t)$,
$1\leq j \leq t$, evaluated at
$0$, form a basis of $\C^t$. Then
\begin{equation}\label{E:neweqn}
\Delta \phi_l^k - \sum_{m=1}^{t} \Delta_{mk}
\phi_l^m \in I, \quad 1\leq l\leq d^\prime, \quad t+1 \leq
k \leq n^\prime,
\end{equation}
where
\begin{equation}\label{E:delta}
\Delta(z,w) =
\begin{vmatrix}
L^{\alpha^1} \phi_{l_1}^1 &\dots &L^{\alpha^1} \phi_{l_1}^t\\
\vdots &\ddots &\vdots\\
L^{\alpha^t} \phi_{l_t}^1 &\dots &L^{\alpha^t} \phi_{l_t}^t
\end{vmatrix}(z,w,0,w), \quad \Delta(0)\neq 0,
\end{equation}
and
\begin{equation}\label{E:deltamk}
\Delta_{mk}(z,w) = (-1)^{t+m}
\begin{vmatrix}
L^{\alpha^1} \phi_{l_1}^1 &\dots &\widehat{ L^{\alpha^1}
\phi_{l_1}^m }&\dots &L^{\alpha^1} \phi_{l_1}^k \\
\vdots & &\vdots & &\vdots \\
L^{\alpha^t} \phi_{l_t}^1 &\dots &\widehat {L^{\alpha^t}
\phi_{l_t}^m} &\dots &L^{\alpha^t} \phi_{l_t}^k
\end{vmatrix}(z,w,0,w),
\end{equation}
and where the $\,\, \widehat{} \,\, $ means that
this  column has been dropped. More specifically,
if $H$ is constantly $(k_0,s)$-degenerate, then
the $\alpha^j$ fulfill $1\leq |\alpha^j| \leq
k_0$, and the same choice for $\alpha^j$ is
possible for every map $H^\prime:M\to M^\prime$
(of constant degeneracy) agreeing with
$H$ up to order $k_0$. 
\end{lem}
\begin{proof}
We will be using the parametrization of $M$ as in
\ref{subsub:normalcoords}. Note that by
\eqref{E:normal1}, for a formal series $\phi
(z,w,\chi,\tau ) \in I$, $\phi (z, w, 0 ,w) = 0$.
We use this in \eqref{E:matrix1}. Developing the
resulting determinant along the last row, we see
that for every $\beta\in\Nn$, for every $l$, $1\leq
l \leq d^\prime$, and for every $k$, $t+1\leq k\leq
n^\prime$
\begin{equation}
\label{E:developedmatrix}
\Delta (z, w) L^\beta \phi_l^k (z,w,0,w) -
\sum_{m=1}^{t} \Delta_{mk}(z,w) L^\beta \phi_l^m
(z,w,0,w)=0,
\end{equation}
where $\Delta$ and $\Delta_{mk}$ are defined
by \eqref{E:delta} and \eqref{E:deltamk},
respectively. Recalling
\eqref{E:dchiasL} we conclude that
\begin{equation}
\label{E:developedmatrix2}
\Delta (z, w)  \phi_l^k (z,w,\chi,\bar Q (\chi, z,w))
-
\sum_{m=1}^{t} \Delta_{mk}(z,w) \phi_l^m
(z,w,\chi,\bar Q (\chi, z,w))=0.
\end{equation}
This immediately implies \eqref{E:neweqn}. The
last statement follows from the construction.
\end{proof}
We are now going to characterize the degeneracy $s$
of a mapping of constant degeneracy in terms of
certain formal holomorphic vector fields. These
results generalize some results about holomorphic
nondegeneracy (defined in \cite{S2}) contained in
e.g. \cite{BERbook}.
\begin{defin} \label{D:holvf} Let $M\subset\CN$,
$M^\prime\subset\CNp$ be formal submanifolds of
$\CN$ and $\CNp$, respectively, $H:M\to M^\prime$ a
formal holomorphic map. A formal holomorphic
vector field $X$ in $\CNp$ tangent to
$M^\prime$ along
$H(M)$ is an operator of the form
\begin{equation}
X = \sum_{j=1}^{N^\prime} a_j (Z)
\vardop{}{Z_j^\prime},
\end{equation}
where $a_j\in\fps{Z}$ (called the
coefficients of $X$), with the property that  for
every
$\phi\in I^\prime$, $\sum_{j=1}^{N^\prime} a_j(Z)
\phi_{Z_j^\prime} (H(Z),\bar H(\zeta)) \in I$.
We say that a set $\{ X_1,\dots
,X_l\}$ of such formal holomorphic vector fields is
linearly independent if their coefficients
evaluated at
$0$ form a linearly independent set of vectors in
$\CNp$.
\end{defin}
This
terminology is motivated by thinking about
$H$ as parametrizing a copy of $\CN$ in $\CNp$, for
which
$Z$ acts as a coordinate function.  In the
special case of an immersion, one can associate
to $X$ a formal holomorphic vector field in
$\CNp$ in the following way. If
$H$ is an immersion, $H(M)$ can be regarded as a
formal submanifold ({\em not} generic) of $\CNp$ by
taking as its ideal $(H^\sharp)^{-1}(I)\supset
I^\prime$. In this case, $H$ has a right inverse
$G$, so
$H^\sharp$ has a left inverse $G^\sharp$. Hence, an
expression of the form
$X=\sum_{j=1}^{N^\prime} a_j(Z) \dop{Z_j^\prime}$
gives rise to a formal holomorphic vector field
$X^\prime=\sum_{j=1}^{N^\prime} G^\sharp
a_j (Z^\prime)
\dop{Z_j^\prime}$ in $\CNp$, and $X$ is tangent to
$M^\prime$ along $H(M)$ if and only if $H^\sharp
(X^\prime f) \in I$ for all $f\in I^\prime$.

Also note that it is enough to check the condition
in Definition~\ref{D:holvf} on a set of generators
of $I$, and that the space of all holomorphic
vector fields tangent to $M^\prime$ along $H(M)$
can be identified with a submodule of
 the free module $\fps{Z}^{N^\prime}$.
\begin{ex}\label{EX:balls} Consider the standard
linear injection of the boundary of the $N$-ball
$|Z_1|^2 +\dots +|Z_N|^2 =1$ into the boundary of the
$N^\prime$-ball $|Z_1|^2
+\dots +|Z_{N^\prime}|^2 =1$,
$(Z_1,\dots,Z_N)\mapsto
(Z_1,\dots,Z_N,0,\dots,0)$. This map is of
constant degeneracy $N^\prime - N$ everywhere and
there are
$N^\prime - N$ linearly independent holomorphic
vector fields tangent to
$M^\prime$ along $H(M)$ given by
$\dop{Z_{N+1}^\prime},\dots,
\dop{Z_{N^\prime}^\prime}$.
\end{ex}
The situation in Example~\ref{EX:balls} is typical
in the following sense:
\begin{prop}\label{P:deginvf}
Assume that $M\subset\CN$, $M^\prime\subset \CNp$
are generic formal submanifolds of $\CN$ and $\CNp$,
respectively, and $H:M\to M^\prime$ is of constant
degeneracy $s$. Then
\begin{multline} \label{E:sisdim}
s=\dim_\C \{X(0)\colon X \text{ formal holomorphic
vector field } \\
 \text{ tangent to } M^\prime \text{ along } H(M)
\}.
\end{multline}
\end{prop}
\begin{proof} Let $\tilde s$ denote the dimension of
the space on the right hand side of
\eqref{E:sisdim}. From Lemma~\ref{L:constdeg} we see
that there are at least $s$ linearly independent
holomorphic vector fields along $H$  tangent to
$M^\prime$ along $H(M)$, so that $s\leq \tilde s$.
On the other hand, assume that $\{ X_1,\dots
X_{\tilde s} \}$ are linearly independent
holomorphic vector fields tangent to $M^\prime$
along $H(M)$; say $X_k=\sum_{j=1}^{N^\prime} a_j^k(Z)
\dop{Z_j^\prime}$, and let $\rho^\prime=
(\rho_1^\prime,\dots , \rho_{d^\prime}^\prime )$ be
a set of generators for $I^\prime$. We have that
\begin{equation}
\label{E:holvfappl}
\sum_{j=1}^{N^\prime} a_j^k (Z) \rho_{l
Z_j^\prime}^\prime (H(Z),\bar H (\zeta )) \in I,
\quad 1\leq l \leq d^\prime.
\end{equation}
Applying CR-vector fields  $L_1,\dots, L_r$
tangent to $M$ to \eqref{E:holvfappl}, we see that
\begin{equation}\label{E:holvfappl2}
\sum_{j=1}^{N^\prime} a_j^k (Z) L_1\cdots L_r \rho_{l
Z_j^\prime}^\prime (H(Z),\bar H (\zeta )) \in I,
\quad 1\leq l \leq d^\prime.
\end{equation}
Evaluating \eqref{E:holvfappl2} at $0$, we conclude
that $\dim_\C E_k (0) \leq N^\prime - \tilde s$ for
all $k$. Hence, $\tilde s \leq s$, and the proof is
complete.
\end{proof}
Note that the second part of the proof of
Proposition~\ref{P:deginvf} shows that if we denote
the dimension of the space on the right hand side of
\eqref{E:sisdim} by $\tilde s$, then the degeneracy
$s$ of a formal map $H$ always satisfies $\tilde s
\leq s$, whether the degeneracy is constant or not.

We now want to relate our notion of nondegeneracy
of a map with the notion of finite nondegeneracy of
manifolds. In particular, we give a bound on the
degeneracy for a certain class of maps between
finitely nondegenerate manifolds.

\subsection{Finitely nondegenerate
manifolds}\label{S:finitenondegeneracy}
The notion of finite nondegeneracy was introduced
in \cite{BHR1}, and has proved very useful. We say
that  a generic submanifold is {\em finitely
nondegenerate} (or, more specifically,
$\ell_0$-nondegenerate) if its identity map is 
$\ell_0$-nondegenerate in the sense of
Definition~\ref{D:formaldegeneracy}. For the
original definition, see e.g. \cite{BERbook},
Chapter IX. By the chain rule we see that if
there is a $k_0$-nondegenerate map into some
generic formal submanifold
$M^\prime\subset\CNp$, then
$M^\prime$ is $\ell_0$-nondegenerate for some
$\ell_0 \leq k_0$. In fact, we also see that every
formal  biholomorphism between generic
formal submanifolds
$M\subset\CN$ and $M^\prime\subset\CN$ of the same
codimension which are $\ell_0$-nondegenerate is
in fact $\ell_0$-nondegenerate.

In order to use finite nondegeneracy of
submanifolds to put bounds on the degeneracy, we
need the mapping to fulfill another crucial
property, which we will introduce next.
\begin{defin}\label{D:transversality}
Let $M\subset\CN$, $M^\prime\subset\CNp$ be formal
generic submanifolds in $\CN$ and $\CNp$,
respectively, $H:M\to M^\prime$ a formal holomorphic
map. We say that $H$ is transversal if for one (and
hence every) set of generators $\rho^\prime =
(\rho_1^\prime ,\dots , \rho_{d^\prime}^\prime)$ of
$I^\prime$, $H^\sharp \rho^\prime$ generates $I$.
\end{defin}
Note that in particular, if $H:M\to M^\prime$ is
transversal, then $d\leq d^\prime$. Transversality
is most easily expressed after choosing normal
coordinates $(z,w)\in\Cn\times\Cd$ for $M$ and
$(z^\prime, w^\prime)\in\Cnp\times\Cdp$ for
$M^\prime$. Write $H=(f,g)$ in these coordinates.
Then $H$ is transversal if and only if the matrix
$\vardop{g}{w}(0)$ has maximal rank $d$.
\begin{lem}\label{L:bound1} Let $M\subset\CN$,
$M^\prime\subset\CNp$ be formal generic submanifolds
of $\CN$ and $\CNp$, respectively, with $M$ being
$\ell_0$-nondegenerate, and $H:M\to M^\prime$ a
transversal mapping. Then the degeneracy $s$ of $H$
fulfills $0\leq s \leq N^\prime - N$.
\end{lem}
\begin{proof} Let $\rho^\prime =
(\rho_1^\prime,\dots,\rho_{d^\prime}^\prime)$
generate $I^\prime$. We can assume that
$(H^\sharp \rho_1^\prime, \dots,
H^\sharp \rho_d^\prime )$ generate $I$. Now using
the chain rule it follows that
\begin{equation}\label{E:rhozis}
\left( \rho_j^\prime (H(Z) , \bar H( \zeta))
\right)_Z = \rho_{j,Z^\prime}^\prime (H(Z),\bar
H(\zeta)) \vardop{H}{Z} (Z), \quad j=1,\dots , d.
\end{equation}
Applying CR-vector fields $L_1,\dots L_r$ tangent
to $M$ to \eqref{E:rhozis}, we see that
\begin{multline}\label{E:rhozis2}
L_1\cdots L_r \left( \rho_j^\prime (H(Z) , \bar H(
\zeta))
\right)_Z = \\ \left( L_1\cdots
L_r\rho_{j,Z^\prime}^\prime (H(Z),\bar
H(\zeta))\right)
\vardop{H}{Z} (Z), \quad j=1,\dots , d.
\end{multline}
By hypothesis, if evaluated at $0$, the dimension of
the space spanned by the vectors on the right hand
side of \eqref{E:rhozis2} is $N$. On the other hand,
the span of the vectors $L_1\cdots
L_r\rho_{j,Z^\prime}^\prime (H(Z),\bar
H(\zeta))$ evaluated at $0$ has dimension $N^\prime
- s$, where $s$ is the degeneracy of $H$. Hence,
\eqref{E:rhozis2} implies that $N\leq N^\prime - s$,
which is the inequality claimed.
\end{proof}
\subsection{Real-analytic and $C^\infty$
Submanifolds}\label{S:raandcinfty} We now want to
apply the theory developed above to concrete
submanifolds of
$\CN$ and $\CNp$. First, let $M\subset\CN$ and
$M^\prime\subset\CNp$ be
generic $C^\infty$-submanifolds of
$\CN$ and $\CNp$ of codimension $d$ and $d^\prime$,
respectively. Assume that $p_0\in M$,
$p_0^\prime\in M^\prime$, and $H$ is a holomorphic
mapping (or, more generally, a
$C^\infty$-CR-mapping) defined in a
neighbourhood $U$ of $p_0$, with $H(p_0) =
p_0^\prime$ and $H(U\cap M)\subset M^\prime$.
We write $\crb{M}$ for the CR-bundle of $M$,
i.e. the bundle with $\crb{M}_p = \C T_p M \cap
\C T^{(0,1)}_p \CN$.

To $M$
and $M^\prime$ we associate formal submanifolds of
$\CN$ and $\CNp$, respectively, by choosing
holomorphic coordinates $Z$ and $Z^\prime$ in $\CN$
and $\CNp$, respectively, in which $p_0 = 0$ and
$p_0^\prime= 0$ and assigning them the ideals
$I\subset \fps{Z,\zeta}$ and $I^\prime \subset
\fps{Z^\prime, \zeta^\prime}$ which are generated
by the Taylor series of their defining
functions. $H$ corresponds to a formal
holomorphic map---by its Taylor expansion, if
it is holomorphic, and by its formal
holomorphic power series (see
\cite{BERbook}, \S 1.7.) if it is 
$C^\infty$-CR. Also, a local basis $L_1,\dots
L_n$ of the CR-vector fields tangent to $M$
gives rise (by taking the Taylor expansion of
the coefficients) to a basis for the formal
CR-vector fields tangent to the formal
manifold M. 

Abusing notation, we shall always use the same
letters to denote the formal object associated
to a concrete object; this will cause no
confusion, since the operations done on them
clearly distinguish the two classes. 

Choose defining
functions
$\rho^\prime = (\rho_1^\prime,\dots
,\rho_{d^\prime}^\prime)$ for 
$M^\prime$ and a local
basis $L_1,\dots , L_n$ for $C^\infty (M,
\crb{M})$. As above, for a multiindex
$\alpha\in\Nn$, we write $L^\alpha =
L_1^{\alpha_1}\dots L_n^{\alpha_n}$.  After
possibly shrinking $U$,  we can define the vector subspaces
\begin{equation}\label{E:theekp}
E_k^\prime (p) = \spanc \{ L^\alpha
\rho^\prime_{j,Z^\prime} (H(Z),\overline{H(Z)} )
\big|_{Z=p} \colon |\alpha |\leq k, 1\leq j \leq
d^\prime \}\subset \CNp
\end{equation}
for $p\in U$. 
Let $s(p)= N^\prime - \max_k \dim_\C E_k^\prime (p)$.
We can then say that $H$ is of degeneracy $s(0)$ at
$0$, and that $H$ is of constant degeneracy $s$ at
$0$ if
$s(p)$ is constant on a neighbourhood of $0$ in $M$.
By taking $k_0$ to be the least integer for which
$E_k^\prime(0) = E_{k_0}^\prime(0)$ for $k\geq k_0$,
we can also define the finer invariant of
$(k_0,s)$-degeneracy, like in
Definition~\ref{D:formaldegeneracy}. Just as in the
case of formal degeneracy, one sees that this
definition is in fact independent of the choices
made, and invariant under biholomorphic changes of
coordinates in both
$\CN$ and $\CNp$. 

Finally, as noted above, the notion of
$k_0$-nondegeneracy makes sense for mappings which
are a priori only assumed to be $C^{k_0}$. This was
used in \cite{Boern2} to prove a reflection
principle, and is used in the statement of
Corollary~\ref{C:jetdetnondeg}.

In the case of real-analytic submanifolds, we can
give generic bounds on both $k_0$ and $s$ (under some
additional assumptions), which we want to do now.

\begin{defin}\label{D:gendeg} Let $M$ and $M^\prime$
be  connected, real-analytic, generic submanifolds
of
$\CN$ and $\CNp$, respectively, and $H$ a
holomorphic mapping defined on an open set $U\in
\CN$ containing $M$ with $H(M)\subset M^\prime$. For
all
$p\in M$, let
$s(p)$ be the degeneracy of $H$ at $p$. The generic
degeneracy $s(H)$ is defined as
$s(H)=
\min_{p\in M} s(p)$. 
\end{defin}
The following Lemma implies that the set of points
where $H$ is of constant degeneracy $s(H)$ is an
open, dense subset of $M$. 
\begin{lem}\label{L:gendeg} Let $M$, $M^\prime$, $H$ be as
in Definition~\ref{D:gendeg}. The set
$\{p\in M\colon s(p) > s(H)
\}$ is real-analytic.
\end{lem}
\begin{proof}
We only prove the the statement locally; a
standard connectedness argument gives the global
statement. After choosing local defining functions,
the points where the degeneracy of $H$ is strictly
bigger than $s(H)$ is given by the vanishing of
determinants with real-analytic entries; see the
arguments before \eqref{E:matrix1}.
\end{proof}
The inequality $0\leq s(H)
\leq N^\prime- d^\prime$ holds trivially. The upper
bound corresponding to Lemma~\ref{L:bound1} is
sharper:
\begin{lem}\label{L:bound2} 
Let $M$, $M^\prime$ and $H$ be as in
Definition~\ref{D:gendeg} and assume in addition
that $M^\prime$ is holomorphically nondegenerate and
that $H$ is transversal. Then $0\leq s(H) \leq
N^\prime - N$. 
\end{lem}
The proof is immediate from Lemma~\ref{L:bound1} and
e.g. Theorem 11.5.1. in \cite{BERbook}. For a bound
on $k_0$, assume for simplicity that $0\in M$, $0\in
M^\prime$, and that $H(M)\subset M^\prime$ with
$H(0) = 0$. Also let
normal coordinates
$(z,w)$ for
$M$ and $(z^\prime, w^\prime)$ for $M^\prime$ with
corresponding real-analytic functions
$Q:\C^{2n+d}\to \C^d$ and
$Q^\prime:\C^{2n^\prime +d^\prime}\to \C^{d^\prime}$
(each defined and convergent in a neighbourhood of
$0\in\C^{2n+d}$ and $0\in\C^{2n^\prime +d^\prime}$)
be chosen. That is, both $Q$ and $Q^\prime$ fulfill
\eqref{E:normal1},  $M$ is given by $w=Q(z,\bar
z,\bar w)$ in a neighbourhood of $0$, and $M^\prime$
is given by $w^\prime = Q^\prime (z^\prime, \bar
z^\prime, \bar w^\prime )$. As in
\ref{sub:formaldegeneracy}, we write $H=(f,g)$ and
set
$Q^\prime_{j z_k^\prime}(f(z,w),\overline{f(z,w)},
\overline{g(z,w)}) = \phi_j^k (z,w,\bar z, \bar w)$.
If $H$ is of constant degeneracy $s$ at $0$, say
$H$ is $(s,k_0)$-degenerate at $0$, then after
reordering we may assume that (writing
$e=n^\prime -s$) the vector valued functions
\begin{equation}\label{E:vectorstospan}
\phi_j = (\phi_j^1,\dots , \phi_j^e, e_j),
\qquad 1\leq j
\leq d^\prime,
\end{equation}
where $e_j$ is the $j$th unit vector in $\Cdp$, are
real-analytic at $0\in\CN$; they are clearly
linearly independent at $0$, and furthermore, if we
choose the basis of CR-vector fields tangent to
$M$
\begin{equation}
L_k = \vardop{}{\bar z_k} +\sum_{j=1}^d \bar
Q_{j\bar z_k} (\bar z, z , w) \vardop{}{\bar w_j},
\quad 1\leq k
\leq d,
\end{equation}
and let $L^\alpha = L_1^{\alpha_1}\cdots
L_n^{\alpha_n}$, then the set $\{L^\alpha \phi_j
\big|_0 \colon \alpha\in\Nn, 1\leq j\leq
d^\prime\} $ spans $\C^{N^\prime - s}$. Now we can
complexify all of these statements. So we let
$\mathcal{M} = \{ (z,w,\chi,\tau) \in U\subset
\Cn\times\Cd\times\Cn\times\Cd \colon \tau_j =
\bar Q_j(\chi,z,w) \}$ where U is a neighbourhood of
$0$ on which $Q$ is convergent be the
complexification of $M$; $\mathcal{M}$ is a
holomorphic submanifold of codimension $d$ in
$\C^{2N}$. We also need the submanifold
$\mathcal{M}_0$ of dimension $n$ defined by
$\mathcal{M}_0= \{ (\chi,\tau)\in\Cn\times\Cd\colon
\tau_j = \bar Q_j (\chi,0,0)\} = \{(\chi,0)\colon
\chi \in \Cn\}$. The complexifications of the $L_k$
are 
\begin{equation}
\mathcal{L}_k=\vardop{}{\chi_k} +\sum_{j=1}^d \bar
Q_{j\chi_k} (\chi, z , w) \vardop{}{\tau_j},
\quad 1\leq k
\leq d,
\end{equation}
and if we denote the complexification of $\phi_j$
again by $\phi_j$, then we have that the set
$\{\mathcal{L}^\alpha
\phi_j
\big|_0 \colon \alpha\in\Nn, 1\leq j\leq
d^\prime\} $, spans $\C^{N^\prime - s}$. Now note
that we can restrict our attention to
$\mathcal{M}_0$ in this statement, since none of the
$\mathcal{L}_k$ differentiates in $z$ or $w$; hence,
we have that
\begin{equation}
\spanc \Big\{ \vardop{^{|\alpha|}}{\chi^\alpha}
\phi_j(0,0,\chi,0) \Big|_{\chi = 0} \colon
\alpha\in\Nn, 1\leq j \leq d^\prime
\Big\} =
\C^{N^\prime -s}.
\end{equation}
We now apply e.g. Lemma 11.5.4. in \cite{BERbook} to
conclude that generically, the derivatives of order
$\leq N^\prime -  s -d^\prime $ suffice; which in
turn implies that generically, $k_0\leq N^\prime -
d^\prime -s$. We summarize:
\begin{lem} Let $M$ and $M^\prime$
be  connected, real-analytic, generic submanifolds
of
$\CN$ and $\CNp$, respectively, and $H$ a
holomorphic mapping defined on an open set $U\in
\CN$ containing $M$ with $H(M)\subset M^\prime$.
Then there exists numbers $s\leq N^\prime
-d^\prime$ and $k_0\leq N^\prime - d^\prime -s$
such that outside some proper real analytic
subvariety of $M$,
$H$ is
$(k_0,s)$-degenerate.  
\end{lem}
\section{Nondegenerate
Mappings}\label{section:nondeg} In this section we
shall discuss nondegenerate mappings. We start with
the ``basic identity'', and in the next
subsection, prove
Theorems~\ref{T:nondegconv},
\ref{T:jetdetnondeg1}, and \ref{T:nondegalg}. 
\subsection{The Basic Identity} We write
$K(t) = |\{
\alpha\in\N^N \colon |\alpha| \leq t \} |$ for the
number of all multiindeces of length less than
$t$.  For a multiindex
$\alpha$, ${\partial^\alpha}$ let  denote the
operator
$\vardop{^{|\alpha|}}{Z^\alpha}$. The following
proposition is our starting point. 
\begin{prop}[Basic Identity for nondegenerate maps]
\label{P:basicidnondeg}
Let $M\subset\CN$, $M^\prime\subset\CNp$ be generic
formal submanifolds, $H:M\to M^\prime$ a formal
holomorphic map which is 
$k_0$-nondegenerate. Then there exists a
formal function $\Psi:
\CN\times\CN\times\C^{K(k_0)N^\prime}\to
\CNp$ (that is, if we write $W$ for the
coordinates in $\C^{K(k_0)N^\prime}$, 
$\Psi\in\fps{Z,\zeta,W}^{N^\prime}$)  with the
property that   
\begin{equation}\label{E:basicidnondeg}
 H(Z) - 
\Psi (Z,\zeta,(\partial^\beta
\bar H (\zeta)- \partial^\beta
\bar H (0) )_{0\leq |\beta |\leq k_0}) \in
I^{N^\prime};
\end{equation}
furthermore, $\Psi$ depends only on
$M$, $M^\prime$, and on the values of
$\partial^\beta H (0)$ for $\beta \leq |k_0|$, such
that if $H^\prime:M\to M^\prime$ is another formal
map with $ \partial^\beta H (0)=\partial^\beta
H^\prime (0)$ for $|\beta| \leq k_0$, then
\eqref{E:basicidnondeg} holds with
$H^\prime$ in place of $H$. If
$M$ and
$M^\prime$ are real-analytic, $\Psi$ is
convergent on a neighbourhood of
the origin.
If $M$ and $M^\prime$ are algebraic, so is
$\Psi$.
\end{prop}
\begin{proof} Choose a basis $L_1,\dots, L_n$ of
the CR-vector fields tangent to $M$ and defining
functions
$\rho^\prime=(\rho_1^\prime,\dots,
\rho_{d^\prime}^\prime)$.
By Lemma~\ref{L:deginbasis}, we can choose
$N^\prime$ multiindeces $\alpha^1,\dots
\alpha^{N^\prime}$ and integers $l^1,\dots,
l^{N^\prime}$ with $0\leq |\alpha^j|\leq k_0$,
$1\leq l^j\leq d^\prime$ for all
$j=1,\dots,N^\prime$ such that 
\begin{equation}
\label{E:nondegjac}
\det \left( 
L^{\alpha^j}\rho_{l^j, Z_k^\prime}^\prime
(H(Z),\bar H (\zeta ))\big|_0
\right)_{\substack{1\leq j \leq N^\prime \\ 
1\leq k \leq N^\prime}} \neq 0.
\end{equation}
We write  
$\Phi_j (Z,\zeta,H(Z), \bar H(\zeta
), (\partial^\beta
\bar H (\zeta ))_{1 \leq |\beta |\leq k_0} ) = 
L^{\alpha^j}\rho_{l^j }^\prime (H(Z),\bar
H (\zeta )) \in I$ for $1\leq j \leq N^\prime$;
using the chain rule, we see that
$\Phi_j \in \fps{Z,\zeta , X,Y }[W]$ where $X=
(X_1,\dots , X_{N^\prime})$ represents the $H_j$,
$Y$ stands for $\bar H$ and $ W $ are variables in
$\C^{({K(k_0)}-1) N^\prime}$. We make a change of
variables by replacing
$W$ by $W + \partial^\beta
\bar H (0) _{0\leq |\beta |\leq k_0}$ and write
again $\Phi_j$ in these new variables; hence, $\Phi_j
 (Z,\zeta,H(Z),\bar H(\zeta), (\partial^\beta \bar
H (\zeta )-\partial^\beta
\bar H (0))_{1 \leq |\beta |\leq k_0} ) \in I$, and
$\Phi_j\in \fps{Z,\zeta,X,Y}[W]$; also, $\Phi_j$
depends only on $M$, $M^\prime$, and on the values
of
$\partial^\beta H (0)$ for $|\beta| \leq k_0$.

Now consider the equations 
\begin{equation}
\label{E:nondegsystem}
\Phi_j (Z,\zeta,X,Y,W) = 0, \quad 1\leq j \leq
N^\prime.
\end{equation}
 We claim that this family of equations has a
unique solution in $X$. In fact, if we compute the
Jacobian of \eqref{E:nondegsystem} with respect to
$X$ at $0$, by the definition of $\Phi_j$ and using
\eqref{E:nondegjac}, we see that the Jacobian
matrix $(\vardop{\Phi_j}{X_k}(0) )_{j,k}$ is
nonsingular. Hence, by the formal implicit
function theorem, there exist unique formal power
series
$\Psi_j\in\fps{Z,\zeta,Y,W}$, $j=1,\dots,
N^\prime$, with 
$\Phi_j ( Z,\zeta,\Psi_1(Z,\zeta,Y,W),\dots
,\Psi_{N^\prime} (Z,\zeta,Y,W), Y, W) = 0$.

 We recall that $\Phi_j
 (Z,\zeta,H(Z),\bar H (\zeta ), (\partial^\beta
\bar H (\zeta )-\partial^\beta
\bar H (0))_{1 \leq |\beta |\leq k_0} ) \in I$; if
we replace $Z$ and $\zeta$ by a parametrization
(as for example in \ref{subsub:normalcoords}) of
$I$, say $Z(x)$ and $\zeta(x)$, we conclude that $\Phi_j
 (Z(x),\zeta(x),H(Z(x)), (\partial^\beta \bar H
(\zeta(x) )-\partial^\beta
\bar H (0))_{1 \leq |\beta |\leq k_0} ) = 0$. It
follows that $H_j(Z(x)) =
\Psi_j(Z(x),\zeta(x),\bar H(\zeta (x)),
(\partial^\beta
\bar H (\zeta(x) )-\partial^\beta
\bar H (0))_{1 \leq |\beta |\leq k_0})$, $1\leq j
\leq N^\prime$. Passing back to the ring
$\fps{Z,\zeta}$, we conclude that $\Psi =
(\Psi_1, \dots, \Psi_{N^\prime})$ fulfills
\eqref{E:basicidnondeg}.  

By construction, the map
$\Phi$ depends only on $M$, $M^\prime$, and
$\partial^\beta H (0)$, $0\leq |\beta |\leq k_0$.
The same choice of $\alpha^1 , \dots ,
\alpha^{N^\prime}$ and $l^1,\dots l^{N^\prime}$
works for every other map $H^\prime$ with
$\partial^\beta H^\prime (0)=\partial^\beta H
(0)$, $|\beta |\leq k_0$. Finally, if $M$ and
$M^\prime$ are real analytic or algebraic, we can
choose the defining functions and the basis of
CR-vector fields to be real-analytic (or
algebraic, respectively) and the last two claims of
Proposition~\ref{P:basicidnondeg} follow since
those classes of maps are closed under application
of the implicit function theorem. 
\end{proof}

 We shall need some formal
vector fields tangent to $M$, which will help us to
exploit \eqref{E:basicidnondeg}. Let $\rho =
(\rho_1,\dots ,\rho_{d} )$ be a
real-analytic defining function for $M$.  After
renumbering, we may assume that 
$\rho_{\Hat \zeta}=\left( \vardop{\rho_j}{\zeta_k}
\right)_{1\leq j,k
\leq d}$ is invertible; set 
\begin{equation}
\label{E:svf}
S_j = \vardop{}{Z_j} - \rho_{Z_j} (\rho_{\Hat
\zeta})^{-1} \vardop{}{\Hat \zeta}, \quad j=1,\dots
,N,
\end{equation}
where $\vardop{}{\Hat \zeta}=
(\vardop{}{\zeta_1},\dots , \vardop{}{\zeta_d})^t$.
Then $S_j$ is a (formal) vector field tangent to
$M$, and its coefficients are convergent, if $M$ is
assumed to be real-analytic, and algebraic
functions if $M$ is assumed to be algebraic. If for
$\alpha\in\N^N$ we write $S^\alpha =
S_1^{\alpha_1}\cdots S_N^{\alpha_N}$, then for
$H\in\fps{Z}$, $S^\alpha H = \partial^\alpha H$.
Applying these vector fields repeatedly to
\eqref{E:basicidnondeg} and using the chain rule we
get the following Corollary to
Proposition~\ref{P:basicidnondeg}.
\begin{cor}\label{C:basicidcor} Under the
assumptions of Proposition~\ref{P:basicidnondeg},
the following holds: For all
$\alpha\in\N^N$, there is a formal function
$\Psi_\alpha: \CN\times\CN\times
\C^{K(k_0+ |\alpha |) N^\prime}\to\CNp$ which is polynomial
in its last $(K(k_0+|\alpha|) - K(k_0))N^\prime$ entries
such that 
\begin{equation}\label{E:basicidnondeg2}
\partial^\alpha H (Z) - \Psi_\alpha ( Z,\zeta,
(\partial^\beta \bar H (\zeta) - \partial^\beta \bar
H (0))_{0\leq |\beta |\leq k_0} , (\partial^\beta
\bar H (\zeta) )_{k_0 < |\beta | \leq k_0 +|\alpha
|})\in I;
\end{equation}
$\Psi_\alpha$ depends only on
$M$, $M^\prime$, and on the values of
$\partial^\beta H (0)$ for $|\beta |\leq k_0$,
such that if $H^\prime:M\to M^\prime$ is another
formal map with $ \partial^\beta H
(0)=\partial^\beta H^\prime (0)$ for $|\beta| \leq
k_0$, then
\eqref{E:basicidnondeg2} holds with
$H^\prime$ in place of $H$. If
$M$ and
$M^\prime$ are real-analytic, $\Psi_\alpha$ is
convergent.
If $M$ and $M^\prime$ are algebraic, so is
$\Psi_\alpha$. 
\end{cor} 
The next step is to repeatedly use
\eqref{E:basicidnondeg2} on the Segre sets. Recall
\eqref{E:segremaps1} and \eqref{E:segremapiszero}.
Hence, choosing normal coordinates $Z=(z,w)$,
$\zeta = (\chi , \tau)$, we have that $f(z,0,0,0) =
0$ for all $f\in I$. Applying this fact to
\eqref{E:basicidnondeg2}, we conclude that 
\begin{equation}\label{E:basiconfirst}
\partial^\alpha H (z,0) = \Psi_\alpha ( z,0,0,0,0 ,
(\partial^\beta
\bar H (0) )_{k_0 < |\beta | \leq k_0 +|\alpha
|}).
\end{equation}
Note that the evaluation occuring causes no
problems, since by Corollary~\ref{C:basicidcor},
$\Psi_\alpha$ is a polynomial with respect to these
variables. Hence the right hand side of
\eqref{E:basiconfirst} defines a formal map
$\Cn\to\CNp$, is convergent if $M$ and $M^\prime$
are real-analytic, and algebraic, if $M$ and
$M^\prime$ are algebraic. This is the case $k=0$ of
the following Corollary (we are using the notation
introduced before
\eqref{E:segremapiszero}):
\begin{cor}\label{C:basiconsegres} For all
$\alpha\in\N^N$, there exists a formal function
$\Upsilon_{k,\alpha}:\C^{kn}\to \CNp$
which depends only on $M$, $M^\prime$, and the
derivatives
$\partial^\beta H (0)$ for $| \beta |\leq
(k+1)k_0+|\alpha |$ such that   
\begin{equation}
\label{E:hypfork} 
\partial^\alpha H(v^{k+1} (z,\xi))
= \Upsilon_{k,\alpha} (z,\xi).
\end{equation}
The dependence on the derivatives is as in
Proposition~\ref{P:basicidnondeg}: If $H^\prime:
M\to M^\prime$ is another formal mapping with 
$\partial^\beta H (0)=\partial^\beta H^\prime (0)$
for
$| \beta |\leq (k+1)k_0 + |\alpha| $, then
\eqref{E:hypfork} holds with $H^\prime$ instead of
$H$.  If
$M$ and
$M^\prime$ are real-analytic,
$\Upsilon_{k, \alpha}$ is convergent on a
neighbourhood of
$0\in\C^{kn}$. If
$M$ and $M^\prime$ are algebraic, so is
$\Upsilon_{k,\alpha}$. 
\end{cor} 
\begin{proof}
We note that \eqref{E:basiconfirst} is just the
case $k=0$. By induction, assume the Corollary
holds for $k<k^\prime$. By
\eqref{E:basicidnondeg2},
\begin{equation}
\begin{split}\label{E:indstep}
\partial^\alpha H
(v^{k^\prime+1}(z^\prime,\xi^\prime)) =
 \Psi_\alpha (&
v^{k^\prime +1}(z^\prime,\xi^\prime),\bar
v^{k^\prime} (\xi^\prime),(\partial^\beta
\bar H (\bar v^{k^\prime} (\xi^\prime)) -
\partial^\beta \bar H (0))_{0\leq |\beta |\leq k_0}
, \\&(\partial^\beta
 (\bar H (v^{k^\prime} (\xi^\prime)) )_{k_0 <
|\beta |
\leq k_0 +|\alpha |}).
\end{split}
\end{equation}
Note that the compositions occuring on the right
hand side are all well defined. We
now plug the induction hypothesis
\eqref{E:hypfork} for
$k = k^\prime - 1$ into \eqref{E:indstep}. In fact,
baring
\eqref{E:hypfork}  and replacing $(z,\xi)$ by
$(\xi^\prime)$, we get that
\begin{equation}
\label{E:barvk} \partial^\beta \bar H (\bar
v^{k^\prime} (\xi^\prime) ) =
\Upsilon_{k^\prime - 1,\beta} (\xi^\prime).
\end{equation}
Now the highest order derivative we need is
$|\beta | = k_0 + |\alpha|$, which by assumption
depends on the derivatives of $H$ of order up to
$(k^\prime+2) k_0 + |\alpha | $, which finishes the
induction. \end{proof}
\subsection{Proof of Theorems \ref{T:nondegconv},
\ref{T:jetdetnondeg1}, and \ref{T:nondegalg}} We
start with Theorem~\ref{T:nondegconv}. We use 
Corollary~\ref{C:basiconsegres} for $k=2k_1-1$,
where $k_1$ is the integer given by
Theorem~\ref{T:bercrit}. Since the
manifolds are assumed to be real-analytic,
$\Upsilon_{2k_1,0}$ is convergent in a
neighbourhood of the origin. By
Theorem~\ref{T:bercrit}, we can choose $(z_0,
\xi_0)$ in this neighbourhood with $v^{2k_1} (z_0,
\xi_0) = 0$ and such that the rank of $v^{2k_1}$ is
$N$ at $(z_0,\xi_0)$. As in the remark after
Theorem~\ref{T:bercrit}, this implies that there is
a holomorphic function $\psi$ defined in a
neighbourhood of $0\in\CN$ such that $\psi(0) =
(z_0,\xi_0)$ and $v^{2k} (\psi (Z)) = Z$. Hence, 
\begin{equation}\label{E:convequal}
H(Z) = H(v^{2k} (\psi(Z)) = \Upsilon_{2k_1 - 1,0}
(\psi (Z)).
\end{equation}
Since the right hand side of \eqref{E:convequal} is
convergent, so is the left hand side. This
completes the proof of Theorem~\ref{T:nondegconv}.

Now assume that $H$ is $C^\infty$-CR and that $M$
and $M^\prime$ are smooth. Its associated formal
holomorphic power series is then a formal
holomorphic map between the formal
submanifolds associated to $M$ and $M^\prime$
(see the remarks in \ref{S:raandcinfty}).  We use
Corollary~\ref{C:basiconsegres} for
$k= k_1 - 1$, where $k_1$ is the integer given by
Theorem~\ref{T:bercrit}. Now set $K=k_1$. Then
Corollary~\ref{C:basiconsegres} implies that 
\begin{equation}
H(v^{k_1}(z,\xi)) = \Upsilon_{k_1 - 1,0} (z,\xi ) =
H^\prime (v^{k_1}(z,\xi)).
\end{equation}
But $\rk (v^k) = N$, which by e.g. Proposition
5.3.5. in \cite{BERbook} implies that $H=
H^\prime$ in the sense of equality of formal power
series, which finishes the proof of
Theorem~\ref{T:jetdetnondeg1}. 

Theorem~\ref{T:nondegalg} follows from
Corollary~\ref{C:basiconsegres} exactly like
Theorem~\ref{T:nondegconv}; we just note that it is
enough to check that $H$ is algebraic on some open
set $U$ containing the point $p_0$ where $H$ is
assumed to be $k_0$-nondegenerate.

\section{Levi-nondegenerate hypersurfaces\\
The  case $N^\prime = N+1$}\label{section:np1}
In this section, we will assume that $M$ and
$M^\prime$ are Levi-nondegenerate (at our
distinguished points). We start with a couple of general
facts.
\subsection{Levi-nondegeneracy}
In normal coordinates, 
which we choose  at our distinguished points $p_0$ and
$p_0^\prime$, $M$ being Levi-nondegenerate means
that we can assume 
 \begin{equation}
Q_{z_j \chi_k}(0,0,0) = \delta_{jk}
\epsilon_k, \quad 1\leq k\leq n,  
\end{equation}
where every $\epsilon_k$ is either $+1$ or  $-1$,
and likewise for $M^\prime$. Here is an easy
technical result about the pullback of the
Leviform by a map $H$ in normal coordinates,
which we will use in the proof of Proposition
\ref{P:np1basicid}.
\begin{lem} \label{L:levitransform}
Let $M\subset\CN$, $M^\prime\subset\CNp$ be
given in normal coordinates by
$w=Q(z,\bar{z},\bar{w})$ and
$w^\prime=Q(z^\prime,\bar{z}^\prime,\bar{w}^\prime)$,
respectively. Assume that
$H=(f,g):(\CN,0)\to (\CNp,0)$ is a formal
holomorphic map, and
$H(M)\subset M^\prime$. Then
\begin{equation} \label{E:levitrans}
g_w(0) Q_{z_j\bar{z}_k}(0) =
\sum_{r,s=1}^{n^\prime} Q^\prime_{z_r^\prime
\bar{z}^\prime_s} (0) f_{r z_j}(0) \overline{f_{s
z_k}(0)},\quad 1\leq j,k\leq n.
\end{equation} 
\end{lem}
To prove this, set in
$g(z,w)=Q^\prime(f(z,w),\bar{f}(\chi,\tau),\bar{g}
(\chi,\tau)$ $\tau=0$, $w=Q(z,\chi,0)$ to obtain
$g(z,Q(z,\chi,0)) = Q^\prime
(f(z,Q(z,\chi,0)),\bar{f}(\chi,0),\bar{g}(\chi,0))$.
Differentiation with respect to $z_j$ and
$\chi_k$ and evaluating at $z=\chi=0$ yields
\eqref{E:levitrans}.  
This has the following consequence:
\begin{cor} Suppose that $M\subset\CN$ and
$M^\prime\subset\CNp$ are (formal) real hypersurfaces which
are Levi-nondegenerate at $p_0$ 
and $p^\prime_0$, respectively,
and $H:(\CN,p_0)\to (\CNp, p_0^\prime)$ is a
(formal) holomorphic map which takes $M$ into
$M^\prime$  and is transversal at $p_0$. Then $H$
is immersive.
\end{cor}
This is easy to see using normal coordinates (which we shall
choose in a way such that $p=0$ and $p^\prime = 0$).
 Since $g_{z^\alpha}(0) = 0$, the differential of $H$ has the following form:
\begin{equation} \partial H (0) = 
\begin{bmatrix} f_{1z_1}(0)  &\dots &f_{1z_n}(0) &f_{1w} (0) \\
                \vdots       &      &\vdots      &\vdots     \\
f_{n^\prime z_1}(0)  &\dots &f_{n^\prime z_n}(0) &f_{n^\prime w} (0) \\
0  &\dots & 0 &g_{w} (0) 
\end{bmatrix}
\end{equation}
and $H$ is immersive if this matrix has rank $N$. Hence, if
$H$ is transversal, $H$ is immersive if and only if the
matrix 
\begin{equation}
\partial f (0) = 
\begin{bmatrix} f_{1z_1}(0)  &\dots &f_{1z_n}(0)  \\
                \vdots       &      &\vdots           \\
f_{n^\prime z_1}(0)  &\dots &f_{n^\prime z_n}(0)  
\end{bmatrix}
\end{equation}
has rank $n$. But by \eqref{E:levitrans}, 
\begin{equation}
g_w (0) \left( Q_{z_j\bar{z}_k}(0) \right)_{1\leq j,k \leq
n} = \overline{\partial f (0)}^t 
\left( Q^\prime_{z_r^\prime
\bar{z}^\prime_s} (0)
\right)_{1\leq r,s \leq n^\prime} \partial f (0).
\end{equation}
This implies that if $H$ is transversal, the rank
of
$\partial f (0)$ is at least $n$, which proves the
corollary.

\subsection{The basic identity for $1$-degenerate
maps} From now on we shall assume that
$N^\prime = N+1$. Note that in the
Levi-nondegenerate case,
\begin{equation} \label{E:lkqzjis} L_k
Q_{z^\prime_j}(f(z,w),\overline{f(z,w)},
\overline{g(z,w)})(0) = \epsilon_j^\prime
\overline{f_{j{z}_k}(0)}, \quad 1\leq k \leq n,\quad 1\leq j
\leq n+1
. \end{equation}
By Lemma
\ref{L:bound1}, if $N^\prime=N+1$ and $H$ is transversal,
the degeneracy $s$ of $H$ at $p_0$ is either $0$ or $1$. 
In the case $s=0$, we can apply Theorem \ref{T:nondegconv} and
\ref{T:jetdetnondeg1} to obtain Theorem \ref{T:nponeconv} and
\ref{T:jetdetnp1}, since by Theorem~\ref{T:bercrit} we see
that $K\leq 2$ if the source manifold is a hypersurface. 
Hence, from now on we will assume that $s=1$. In this
subsection, we will develop a basic identity for
1-degenerate maps between Levi-nondegenerate
hypersurfaces.
From \eqref{E:lkqzjis} we see that in Lemma~\ref{L:constdeg}
we can choose $\alpha^j$ to be the multiindex with a
$1$ in the $i$-th spot and $0$ elsewhere and
reorder the
$z^\prime$'s, to get that after barring
\eqref{E:neweqn},
\begin{equation}\label{E:neweqnbar} \begin{split}
\bar \Delta (\chi, \tau)& \bar
Q^\prime_{\chi_{n+1}^\prime} (\bar f (\chi,\tau) ,
f(z,w) , g(z,w)) )  -  \\  & \sum_{m=1}^n \bar
\Delta_m (\chi, \tau) \bar
Q^\prime_{\chi_{m}^\prime} (\bar f (\chi,\tau) ,
f(z,w) , g(z,w)) ) \in I,
\end{split}
\end{equation} 
where now
\begin {equation} \label{E:deltafornp1}
\Delta (z,w) = 
\begin{vmatrix}
L_1 Q^\prime_{z_1^\prime} (f,\bar f, \bar g)& \dots
&L_1 Q^\prime_{z_n^\prime} (f,\bar f, \bar g)  \\
\hdots & & \hdots \\
L_n Q^\prime_{z_1^\prime} (f,\bar f, \bar g)& \dots
&L_n Q^\prime_{z_n^\prime} (f,\bar f, \bar g) 
\end{vmatrix} (z,w,0,w) ,
\end{equation}
\begin{equation} \label{E:deltafornp1at0}
\quad \bar \Delta (0) = 
\begin{vmatrix} \epsilon_1^\prime  f_{1z_1}(0)  &\dots
&\epsilon_n^\prime f_{n z_1}(0)  \\                 \vdots       &      &\vdots           \\
\epsilon_1^\prime f_{1 z_n}(0)  &\dots &\epsilon_n^\prime f_{n z_n}(0)  
\end{vmatrix} \neq 0,
\end{equation}
\begin{equation} \label{E:deltamfornp1at0}
\bar \Delta_m (0) = (-1)^{n+m}
\begin{vmatrix} 
\epsilon_1^\prime  f_{1z_1}(0)  &\dots
&\widehat{\epsilon_m^\prime f_m z_1 (0)} &\dots
&\epsilon_{n^\prime}^\prime f_{n^\prime z_1}(0)  \\         
\vdots       &      &\vdots   & &\vdots        \\
\epsilon_1^\prime f_{1 z_n}(0)  &\dots
&\widehat{\epsilon_m^\prime f_m z_n (0)} &\dots
&\epsilon_{n^\prime}^\prime f_{n^\prime z_n}(0)  
\end{vmatrix}.  \end{equation}
Since $H$ maps $M$ into $M^\prime$,  the chain rule implies
that we 
have formal functions $\Phi_1,\dots , \Phi_n$ such that 
$\Phi_j (Z,\zeta, H(Z), \bar H(\zeta) 
,\partial \bar H (\zeta)) = L_j Q^\prime (f(z,w), 
\bar f (\chi,\tau) , \bar
g(\chi,\tau))$ which are convergent if $M$ and $M^\prime$
are, and are polynomial in the derivatives of $\bar H$. As
in the proof of Proposition~\ref{P:basicidnondeg} we obtain
functions $\Phi_j\in \fps{Z,\zeta,X,S,W}$ where $X,S \in
\C^{N+1} $, $W\in \C^{(N+1)^2}$, such that
\begin{equation} \label{E:phijis}
\Phi_j (Z,\zeta, H(Z), \bar H(\zeta) ,\partial \bar H
(\zeta) - \partial \bar H (0) ) \in I.
\end{equation}
From \eqref{E:neweqnbar} we conclude that there is a formal
function $\Upsilon\in \fps{Z,\zeta,X,S,T,W^\prime}$
such that 
\begin{equation} \label{E:theupsilon}
\Upsilon (Z,\zeta, H(Z), \bar H( \zeta), H(\zeta) ,
\partial H (0,\tau)
 - \partial H (0)) \in I.
\end{equation} 
$\Phi_1,\dots, \Phi_n , \Upsilon$ only depend on $M$,
$M^\prime$ and the derivative of $H$ evaluated at $0$,
will agree for $H$ and $H^\prime$ with $\partial H(0) =
\partial H^\prime (0)$, are convergent if $M$ and $M^\prime$
are real-analytic, and algebraic if $M$ and
$M^\prime$ are algebraic. 
 Consider the system of equations  
 \begin{equation}\label{E:systemfornp1}
 \begin{split} &\Phi_j
 (Z,\zeta,X,S,W  )  =  0, \quad 1\leq j \leq n, \\ 
 & \Upsilon
 (Z,\zeta,X,S,T,W^\prime  ) = 0, \\
 &X_{N+1} = Q^\prime (X^\prime, S),
\end{split}
\end{equation}
in $\fps{Z,\zeta,X,S,T,W,W^\prime }$, where
$X=(X^\prime, X_{N+1})$ is the splitting
corresponding to $H=(f,g)$.  We claim that we can
apply the implicit function theorem to 
\eqref{E:systemfornp1} to see that this system admits a
unique solution in $X$.  
In order to compute the Jacobian of \eqref{E:systemfornp1}
with respect to $X$, first note that since
$\Upsilon_{X_{N+1}
} (0,0,0,0,0) = 0$ and $X_{N+1}$ does not appear in any of the
$\Phi_j$, it is enough to show that the determinant
\begin{equation} D = 
\begin{vmatrix}
\Phi_{1 X_1}(0) &\dots &\Phi_{1 X_{N}} (0) \\
\vdots &                &\vdots \\
\Phi_{n X_1}(0) &\dots &\Phi_{n X_{N}} (0) \\
\Upsilon_{X_1} (0) &\dots &\Upsilon_{X_N} (0)
\end{vmatrix}
\end{equation}
is nonzero. Note that $\Phi_{ k X_j} = \epsilon_j^\prime
\overline{f_{j z_k} (0)}$, that $\Upsilon_{X_j}(0) = -
\epsilon_j^\prime \bar \Delta_j (0)$ for $1\leq j \leq n$
and $\Upsilon_{X_N} (0) = \epsilon_N^\prime \bar \Delta
(0)$. To simplify notation in the following
argument, we write
$\Delta (0) = -
\Delta_{N} (0)$ for ease of notation. Developing
$D$ along the last row and using
\eqref{E:deltafornp1at0} and \eqref{E:deltamfornp1at0}
we see that  \begin{multline} \label{E:developeddet}
D = \pm \left(
\epsilon_1^\prime 
\begin{vmatrix}
     	f_{2 z_1} &\dots &f_{n^\prime z_1}\\
      \vdots    &      &\vdots          \\
      f_{2 z_n} &\dots &f_{n^\prime z_n} 
\end{vmatrix}
\begin{vmatrix} 
      \overline{f_{2 z_1}} &\dots &\overline{f_{n^\prime z_1}}\,\\
      \vdots               &      &\vdots                    \\
      \overline{f_{2 z_n}} &\dots &\overline{f_{n^\prime z_n}}\,
\end{vmatrix} \right.
\\ +\epsilon_2^\prime 
\begin{vmatrix} 
     f_{1 z_1}  &f_{3 z_1} &\dots &f_{n^\prime z_1} \\
     \vdots     &\vdots    &      &\vdots           \\
     f_{1 z_n}  &f_{3 z_n} &\dots &f_{n^\prime z_n} 
\end{vmatrix}
\begin{vmatrix} 
     \overline{f_{1 z_1}} &\overline{f_{3 z_1}} &\dots &\overline{f_{n^\prime z_1}}\,\\
     \vdots               &\vdots               &      &\vdots                     \\
     \overline{f_{1 z_n}} &\overline{f_{3 z_n}} &\dots &\overline{f_{n^\prime z_n}}\,
\end{vmatrix}
+\dots 
\\ \dots
+\epsilon_{n^\prime}^\prime \left.
\begin{vmatrix}
     f_{1 z_1}  &\dots  &f_{n-1 z_1} &f_{n z_1} \\
     \vdots     & &\vdots            &\vdots    \\
     f_{1 z_n}  &\dots  &f_{n-1 z_n} &f_{n z_n}
\end{vmatrix} 
\begin{vmatrix} 
     \overline{f_{1 z_1}} &\dots &\overline{f_{n-1 z_1}} &\overline{f_{n z_1}}\,\\
     \vdots               &      &\vdots                 &\vdots              \\
     \overline{f_{1 z_n}} &\dots &\overline{f_{n-1 z_n}} &\overline{f_{n z_n}}\,
\end{vmatrix} \right)
\end{multline}
We apply the  Cauchy-Binet Formula to
\eqref{E:developeddet} to see that $\pm D$ is
equal to the determinant of
\begin{equation}
\label{E:prod} \begin{bmatrix} 
     f_{1 z_1}  &f_{2 z_1} &\dots &f_{n^\prime z_1} \\
     \vdots     &\vdots    &      &\vdots           \\
     f_{1 z_n}  &f_{2 z_n} &\dots &f_{n^\prime z_n} 
\end{bmatrix}
\begin{bmatrix}
					\epsilon_1^\prime  &0                 &\dots &0 \\
					0                  &\epsilon_2^\prime &\dots &0 \\
     \vdots             &\vdots            &\ddots &\vdots \\
     0                  &0                  &\dots  &\epsilon_{n^\prime}^\prime
\end{bmatrix}
\begin{bmatrix}
     \overline{f_{1 z_1}} &\dots  &\overline{f_{1 z_n}} \\
     \overline{f_{2 z_1}} &\dots  &\overline{f_{2 z_n}} \\
					\vdots               &        &\vdots               \\
     \overline{f_{n^\prime z_1}} &\dots  &\overline{f_{n^\prime z_n}} 
\end{bmatrix}
\end{equation}
Now apply Lemma \ref{L:levitransform} 
to see that the determinant of 
\eqref{E:prod} is just $\pm g_w (0)$
which we assume to be nonzero. Hence, the claim is proved,
and summarizing, we have proved the following.

\begin{prop} \label{P:np1basicid}
Let $M\subset\CN$ and $M^\prime\subset\C^{N+1}$ be
Levi-nondegenerate formal real hypersurfaces. Let
$H: M \to M^\prime$ be a formal holomorphic map
which is constantly $1$-degenerate and
transversal. Let
$Z=(z,w)$,
$\zeta = (\chi,\tau)$ be normal coordinates for $M$. Then
there exists a formal function $\Psi: \CN\times\CN\times
\C^{N+1}\times \C^{N+1} \times \C^{N(N+1)} \times
\C^{N(N+1)}
\to \C^{N+1}$ such that  \begin{equation}
\label{E:basicidnp1} H(Z) - \Psi (Z,\zeta,\bar H
(\zeta), H (\zeta),
\partial
\bar H(\zeta) -\partial \bar H (0),
\partial H(0,\tau) - \partial H (0) ) \in I^{N+1}.
\end{equation}
Furthermore, $\Psi$ depends only on $M$, $M^\prime$ and the
first derivative of $H$ at $0$, such that if $H^\prime$ is
another map fulfilling the assumptions of the proposition
with $\partial H (0) = \partial H^\prime (0)$ then
\eqref{E:basicidnp1} holds with $H$ replaced by $H^\prime$. 
If $M$ and $M^\prime$ are real-analytic, $\Psi$ is
convergent on a neighbourhood of the origin. If $M$ and
$M^\prime$ are algebraic, so is $\Psi$.
\end{prop} 
Differentiating this identity as in the proof of
Corollary~\ref{C:basicidcor} we obtain the following.
\begin{cor} \label{C:np1basicidcor} Under the assumptions of
		Proposition~\ref{P:np1basicid}, the following holds:
		For all multiindeces $\alpha\in\Nn$, 
		there 		is 		a 		formal 		function 		
	$\Psi_\alpha:\CN\times\CN\times\C^{N+1}\times
\C^{N+1}\times \C^{(K(1+|\alpha |)-1)
		(N+1)} \times\C^{K(|\alpha|)}
\times \C^{(K(1+|\alpha |)-1)
		(N+1)}\to
		\C^{N+1}$ which is polynomial in its $6$th,
$7$th and $9$th variable such that
\begin{equation}\label{E:basicidalphanp1}
\begin{split}
\partial^\alpha H(Z) - 
&\Psi_\alpha (Z,\zeta, \bar H(\zeta),H(\zeta ),
\partial 
\bar H (\zeta ) - \partial\bar H(0), (\partial^\beta \bar
H(\zeta))_{2\leq |\beta | \leq 1+ |\alpha|},\\ 
&(\partial^\beta 
H(\zeta))_{1\leq |\beta | \leq |\alpha|},
\partial 
 H (0,\tau ) - \partial H(0), (\partial^\beta 
H(0,\tau))_{2\leq |\beta | \leq 1+ |\alpha|}) \in I^{N+1};
\end{split}
\end{equation}
$\Psi_\alpha$ depends only on $M$, $M^\prime$, and on the
first derivative of $H$ at $0$ such that if $H^\prime: M\to
M^\prime$ 
is another formal map fulfilling the assumptions of
Proposition~\ref{P:np1basicid} with $\partial H (0) =
\partial H^\prime (0)$ then \eqref{E:basicidalphanp1} holds
with $H$ replaced by $H^\prime$. If $M$ and $M^\prime$ are
real analytic, $\Psi_\alpha$ is convergent. If $M$ and
$M^\prime$ are algebraic, so is $\Psi_\alpha$.
   \end{cor}
The main difference between
\eqref{E:basicidnondeg2} and \eqref{E:basicidalphanp1} is
that in \eqref{E:basicidalphanp1} the argument $(0,\tau)$
appears. This means that we can only iterate
\eqref{E:basicidnp1} once, and hence we can determine $H$
from its $2$-jet at $0$ {\em only on the $2$nd
Segre set}. This is the main reason why we have to
restrict to  hypersurfaces here.
\begin{cor}\label{C:basicidcornp1} Under the assumptions of
Proposition \ref{P:np1basicid}, for all $\alpha\in\N^N$ there is a
 formal function $\Upsilon_\alpha:\Cn\times\Cn\to\C^{N+1}$
 such that
 \begin{equation}\label{E:basicidcornp1}
 \partial^\alpha H (z, Q(z,\chi,0)) = \Upsilon_\alpha
 (z,\chi).
 \end{equation}
$\Upsilon_\alpha$ depends only on $M$, $M^\prime$ and
$\partial^\beta 
H(0)$ for $|\beta | \leq 2+|\alpha |$ such that if
$H^\prime$ is another map fulfilling the hypotheses of
Proposition \ref{P:np1basicid} with $\partial^\beta 
H(0)=\partial^\beta 
H^\prime(0)$ for $|\beta | \leq 2+|\alpha |$, then
\eqref{E:basicidcornp1} holds with $H$ replaced by $H^\prime$. If $M$ and $M^\prime$
are real-analytic, then $\Upsilon_\alpha$ is
convergent on a neighbourhood of
$0\in\Cn\times\Cn$. If $M$ and $M^\prime$ are
algebraic, so is $\Upsilon_\alpha$.
 \end{cor}
The proof of this corollary is by induction just as in
Corollary \ref{C:basicidcor} and left to the reader.
Theorem \ref{T:jetdetnp1} in the case $s=1$ follows
from Corollary \ref{C:basicidcornp1} just as Theorem
\ref{T:jetdetnondeg1} follows from Corollary
\ref{C:basicidcor}. Theorem \ref{T:np1alg} also follows
easily from Corollary~\ref{C:basicidcornp1} since by Lemma
\ref{L:gendeg} we can always pass to a point where $H$ is of
constant degeneracy.
However, since we can only work on the second Segre set, we have to work a
little harder for Theorem \ref{T:nponeconv}. We
are basically following the argument given in
\cite{BER3}.
\subsection{Proof of Theorem \ref{T:nponeconv}} From
Corollary \ref{C:basicidcornp1} we conclude that 
\begin{equation}\label{E:Hisequal}
H(z,w) = \sum_{j=1}^\infty H^j (z) w^j
\end{equation}
where $H^j (z) = \frac{1}{j!} H_{w^j} (z,0)$ is convergent.
We now want to solve the equation $w = Q(z,\chi, 0)$ for
$\chi$ in $w$. Choose a $\chi_0$ such that the function
$\phi (z,t) = Q(z,\chi_0 t, 0)$, which is defined
in a neighbourhood of the origin in $\Cn\times\C$, 
has a derivative in $t$ which is not constantly vanishing. 
We write
\begin{equation}\label{E:phiisequal}
\phi (z,t) = \sum_{j=1}^\infty \phi_j (z) t^j
 \end{equation}
 and define a convergent power series $\psi( z, t) =
 t+\sum_{j=2}^\infty C_j (z) t^j$ where $C_j(z) = \phi_j (z)
 \phi_1(z)^{j-2}$ for $j\geq 2$. By the implicit function
 theorem, the equation $w = \psi(z,t)$ has a
solution $t = \theta (z,w)$ which is convergent in
a
 neighbourhood of the origin in $\Cn\times\C$.
Then $t = \phi_1 (z) \theta (z, \frac{w}{\phi_1(z)^2}) $
solves $w= Q(z,\chi_0 t, 0 )$. By changing $\theta$, we
can assume that $\phi_1 (z) = A(z)$ is a Weierstrass
polynomial. Hence we conclude that 
\begin{equation} \label{E:Hequalconv} H(z,w) =
F(z,\frac{w}{A(z)^2})  \end{equation}
where $F$ is now a function which is convergent in a
neighbourhood of the origin in $\CN$. We expand $F$ in the
following way: $F(z,t) = \sum_{j=1}^{\infty} F_j(z) t^j$.
Comparing coefficients in \eqref{E:Hequalconv}, we conclude
that $H_j(z) = F_j(z) A(z)^{-2j} $. We now apply the
division theorem to see that
\begin{equation}\label{E:Fjdivided}
F_j(z) = B_j (z) A(z)^{2j} + r_j(z)
\end{equation}
where $r_j (z)$ is a ($\C^{N+1}$-valued)
Weierstrass polynomial of degree less than $2jp$
where $p$ is the degree of the Weierstrass
polynomial $A(z)$. Furthermore, we have the
inequality 
\begin{equation}\label{E:normbjest}
\|
B_j(z) \| \leq C^j \|F_j(z) \|
\end{equation}
which holds for $z$ in a neighbourhood of the origin, with
some constant $C$ (see e.g. \cite{Ho1}, Theorem
6.1.1.). Since
$H_j$ is convergent, we conclude that $r_j$ is the
zero polynomial. So
$H_j (z)= B_j(z)$ and from \eqref{E:normbjest} we
finally conclude that $H(z,w)$ is convergent in a
neighbourhood of the origin.

\section{Strictly pseudoconvex
targets}\label{section:pcx} We will just
indicate how to derive a basic identity in this
case; the  proof is then finished by exactly the
same arguments as  in the Levi-nondegenerate case.
By the Chern-Moser normal form (\cite{CM}), we can
in particular  assume that the target hypersurface
is given in normal coordinates
$(z^\prime , w^\prime)$ by $w^\prime = Q^\prime 
(z^\prime, \chi^\prime, \tau^\prime)$, where
\begin{equation}
\label{E:cmassume}
Q^\prime (z^\prime, \chi^\prime , \tau^\prime ) =
\innprod{z^\prime}{\chi^\prime} + \sum_{\substack{
\alpha, \beta, \gamma \\
|\alpha|, |\beta| \geq 2}}
c_{\alpha,\beta,\gamma} {z^\prime}^\alpha
{\chi^\prime}^\beta {\tau^\prime}^\gamma.
\end{equation}
In this equation,  $\innprod{z^\prime}{\chi^\prime}
=
\sum_{j=1}^{n^\prime} z^\prime_j
\chi^\prime_j =
(z^\prime)^t
\chi^\prime$. It follows that 
\begin{equation}\label{E:Lalphaisfbar}
L^\alpha
Q^\prime_{Z^\prime} (f(z,w), \bar f(\chi,\tau),
\bar g(\chi,\tau) ) \big|_{0} = L^\alpha \bar
f(0).
\end{equation}
Hence, if $H$ is constantly $(k_0,s)$
degenerate  at
$0$, we can choose $t = n^\prime - s$ multiindeces
$\alpha^1,\dots, \alpha^t$, $\alpha^j \in \Nn$,
$|\alpha^j|\leq k_0$, such that the vectors
$\xi_j = L^{\alpha^j}\bar f (0)$, $1\leq j\leq t$
are linearly independent in $\Cnp$. We extend this
set to a basis $\xi_1,\dots,\xi_{n^\prime}$ of
$\Cnp$. By the Gram-Schmidt orthogonalization
process, we obtain 
vectors $v_1, \dots, v_{n^\prime}$ which are
orthonormal with respect to the standard hermitian
product on
$\Cnp$, and a lower triangular invertible matrix
$C$ such that $V=CE$, where $V$ denotes the
unitary matrix with rows $v_1, \dots,
v_{n^\prime}$, and $E$ denotes the matrix with
rows $\xi_1,\dots,
\xi_{n^\prime}$. We  change  coordinates by
$\tilde z = Vz$, $\tilde w = w$. Since $V$ is
unitary, the defining function still has the form
\eqref{E:cmassume}; in particular, 
\eqref{E:Lalphaisfbar} holds with $f$ replaced by
$\tilde f$, and $L^\alpha \overline{{\tilde f} (0)}
=
\bar V L^\alpha \bar f(0)$. Since $\bar V E =
(C^t)^-1$ is upper triangular, it follows 
that we can assume   
$\xi_j^k = (L^{\alpha^j}\bar f (0))_k = 0$ for
$j>k$. Note that the same change of coordinates
works for every other map whose $k_0$-jet at the
origin agrees with the $k_0$-jet of $H$, if $k_0
>1$. 

We now start as in the proof of
Proposition~\ref{P:basicidnondeg} and obtain formal
functions
$\Phi_j (Z,\zeta,  X,Y,W )\in \fps{Z,\zeta,
X,Y,W}$, defined by $\Phi_j (Z,\zeta, H(Z) ,\bar H
(\zeta ), (\partial^\beta \bar H(\zeta) - 
\partial^\beta \bar H(0))_{|\beta |\leq k_0} ) =
(L^{\alpha^j} Q^\prime (f,\bar f, \bar g))
(Z,\zeta) $, which are convergent if $M$ and
$M^\prime$ are real-analytic and algebraic if $M$
and $M^\prime$ are algebraic. Also note that
$\Phi_j$ does not depend on $X_{N^\prime}$. 

The missing equations are supplied by
Lemma~\ref{L:constdeg}. So we define formal
functions $\Upsilon_k (Z,\zeta,X,Y,S,T,W^\prime
)$, $1\leq k \leq s$, by 
\begin{multline}\label{E:theupsilons}
\Upsilon_{j-t} (Z,\zeta, H(Z), \bar H(\zeta),
H(\zeta ), H (0,\tau),(\partial^\beta H (0,\tau) -
\partial^\beta H (0))_{1\leq |\beta|\leq k_0}) =\\
\bar \Delta (\chi,\tau ) \bar Q_{\chi^\prime_{j}}
(\bar f(\chi,\tau ) , f(z,w), g(z,w))  -
\sum_{m=1}^{t}
\bar \Delta_{m} (\chi,\tau )
\bar Q_{\chi^\prime_{m}}
(\bar f(\chi,\tau ) , f(z,w), g(z,w))
\end{multline} 
for $t+1\leq j \leq n^\prime$, and where we let
$\bar\Delta$ be defined by \eqref{E:delta} and
$\bar\Delta_m$ be defined by \eqref{E:deltamk}.
Recall that the functions $\Upsilon_k$ are
convergent if $M$ and $M^\prime$ are
real-analytic and algebraic if $M$ and $M^\prime$
are algebraic. We now claim that we can apply the
implicit function theorem to solve the system
\begin{align}\label{E:pcxsystem}
&\Phi_j (Z,\zeta,X,Y,W) = 0, \quad 1\leq j\leq t,
\notag \\ &\Upsilon_k  (Z,\zeta,X,Y,S,T,W^\prime)
= 0 ,
\quad 1\leq k\leq s,\\
&X_{N^\prime} = Q^\prime( X^\prime,S), \notag
\end{align}
uniquely in $X=(X^\prime, X_{N^\prime})$. First
note that
$\Phi_{j X_{N^\prime}} (0) = \Upsilon_{k
X_{N^\prime}} (0) = 0$ for $1\leq j\leq t$, and
$1\leq k \leq s$. So we only need to  consider the
Jacobian of 
$(\Phi_1,\dots,\Phi_t,\Upsilon_1,\dots,\Upsilon_s)$
with respect to $X^\prime = (X_1,\dots,
X_{n^\prime})$. Now $\Phi_{j X^\prime} (0)  =
\xi_j$, and $\Upsilon_{kX^\prime} (0) =
(\bar\Delta_1 (0),\dots, \bar\Delta_t (0),0,\dots
\bar\Delta (0) , \dots,0)$, where the $\Delta(0)$
appears in the $(t+k)$-th spot, $1\leq k\leq
s$. Recall that $\xi_j^k= 0$, $j>k$, by our
choice of coordinates. Writing out the
determinant we see that indeed the implicit
function theorem applies. This gives the
desired basic identity.    
\bibliographystyle{plain}
\bibliography{deg}
\end{document}